\documentclass[fleqn,12pt]{article}

\usepackage{subcaption,graphicx}
\usepackage{color}
\usepackage{amsfonts, amsmath, amsthm, amssymb,epstopdf}
\usepackage{comment}
\usepackage{sectsty}
\usepackage[margin=2cm]{geometry}
\usepackage[round]{natbib}
\usepackage{setspace}
\usepackage{placeins}
\usepackage[titletoc, toc,page]{appendix}

\usepackage{xr}
\usepackage{bm}

\pdfoutput=1

\usepackage{color}
\definecolor{Blue}{rgb}{0,0,1}
\definecolor{DarkBlue}{rgb}{0,0,0.5}
\definecolor{Red}{rgb}{1,0,0}
\definecolor{Green}{rgb}{0,1,0}
\definecolor{Yellow}{rgb}{1,1,0}
\definecolor{DarkGreen}{rgb}{0,0.4,0}
\definecolor{DarkRed}{rgb}{0.5,0,0}
\definecolor{DarkYellow}{rgb}{0.7,0.7,0}

\usepackage{hyperref}
\hypersetup{
    unicode=false,          					
    pdftoolbar=true,       						
    pdfmenubar=true,        					
    pdffitwindow=true,     						
    pdfstartview={FitH},    					
    pdftitle={Inference on a New Class of Sample Average Treatment Effects}, 
    pdfauthor={Juliana Londono},     		
    pdfsubject={LaTeX},   						
    pdfproducer={Juliana Londono}, 			
    pdfkeywords={}, 
    pdfdisplaydoctitle=true,					
    pdfnewwindow=true,      					
    colorlinks=true,       						
    linkcolor=DarkBlue,          			
    citecolor=DarkBlue,	        			
    filecolor=DarkBlue,      					
    urlcolor=blue,           					
    citebordercolor=0 0 1,	 					
    linkbordercolor=0 0 1,						
    urlbordercolor=0 0 1,							
		frenchlinks=false,    						
		linktocpage=true,
    draft=false
}

\usepackage[para,online,flushleft]{threeparttable}

\graphicspath{{/home/yotam/Desktop/Dropbox/att_CI/figures/}{../../figures/}}

\usepackage{bbm}

\title{ \Large \textbf{
Inference on a New Class \\ of Sample Average Treatment Effects}\thanks{Jasjeet S. Sekhon: 
Robson Professor of Political Science and Statistics, University of California Berkeley; 
\url{sekhon@berkeley.edu};
Yotam Shem-Tov: Ph.D.\ Candidate,
Department  of  Economics,
\url{shemtov@berkeley.edu}.
A previous version of this paper was entitled ``Efficient Estimation of Average Treatment Effects under Effect Heterogeneity.'' We thank Peter Aronow, Max Balandat, Eytan Bakshy, Avi Feller, Johann Gagnon-Bartsch, Peng Ding, Ben Hansen, Guido Imbens, James Robins, S\"{o}ren Kunzel, Winston Lin, Juliana Londo\~{n}o V\'{e}lez, Fredrik S\"{a}vje, and John Myles White for helpful comments and discussions. In addition, we thank Max Balandat, Eytan Bakshy, and John Myles White for help with one of the applications. We also thank the participants of the Berkeley Statistics Annual Research Symposium 2017 and the Atlantic Causal Inference Conference 2017.
Sekhon wishes to acknowledge Office of Naval Research (ONR) grants N00014-15-1-2367 and N00014-17-1-2176.
The R package, \texttt{estCI}, that implements the estimation methods described in the paper is available at \url{https://github.com/yotamshemtov/estCI}. 
}
} 

\author{ Jasjeet S. Sekhon  \\ (UC Berkeley) \and 
Yotam Shem-Tov  \\  (UC Berkeley) }

\newcommand\independent{\protect\mathpalette{\protect\independenT}{\perp}}
\def\independenT#1#2{\mathrel{\rlap{$#1#2$}\mkern2mu{#1#2}}}
 
\theoremstyle{definition}

\newcommand{\e}[1]{\ensuremath{\mathbb{E} \left[ #1 \right]}} 
\newcommand{\ee}[2]{ \ensuremath{  \mathbb{E}_{\text{ #2 }} \left[  #1 \right]  } }
\newcommand{\var}[1]{\ensuremath{\text{Var} \left( #1 \right)}}
\newcommand{\vv}[2]{\ensuremath{\text{Var}_{\text{#2}} \left( #1 \right)}}
 
\newcommand{\cov}[1]{\ensuremath{\text{Cov} \left( #1 \right)}}

\newcommand{\diff}[1]{\ensuremath{ t_{\text{diff}}} \left( #1 \right)}
\newcommand{\td}{\ensuremath{ t_{\text{diff}} } }

\newcommand{\att}{\text{SATT}}
\newcommand{\atc}{\text{SATC}}
\newcommand{\ate}{\text{SATE}}
\newcommand{\ato}{\text{SATO}}

\newtheorem{theorem}{Theorem}[section]

\newtheorem{lemma}{Lemma}

\begin{document}


\maketitle


\doublespacing

\begin{abstract}

We derive new variance formulas for inference on a general class of estimands of causal average treatment effects in a Randomized Control Trial (RCT). We demonstrate the applicability of the new theoretical results using an empirical application with hundreds of online experiments with an average sample size of approximately one hundred million observations per experiment. We generalize \cite{robins1988} and show that when the estimand of interest is the Sample Average Treatment Effect of the Treated (SATT, or SATC for controls), a consistent variance estimator exists. Although these estimands are equal to the Sample Average Treatment Effect (SATE) in expectation, potentially large differences in both accuracy and coverage can occur by the change of estimand, even asymptotically. Inference on SATE, even using a conservative confidence interval, provides incorrect coverage of SATT or SATC. We derive the variance and limiting distribution of a new and general class of estimands---any mixing between SATT and SATC---of which SATE is a specific case. An R package, \texttt{estCI}, that implements all the proposed estimation procedures is available. 


\end{abstract}


\vspace{10 mm}


\section{Introduction}

The Neyman variance estimator is the most commonly used variance estimator in randomized experiments \citep{imbens2015}. Under the super-population model, it is a consistent estimator for the variance of the difference-in-means, and this probably accounts for its popularity. However, under Neyman's finite population model, a consistent variance estimator for the difference-in-means does not exist \citep{neyman1923}, and Neyman's variance estimator is conservative. Sharper, albeit still conservative, variance estimators exist \citep{aronow2014sharp}, but they are not often used.    

To estimate the Sample Average Treatment Effect (SATE), researchers use the limiting distribution of the difference-in-means recentered around SATE. We show that changing the estimand and recentering the difference-in-means limiting distribution with respect to the Sample Average Treatment Effect for the Treated (SATT, or SATC for controls) allows one to obtain a consistent non-conservative variance estimator. 
Consequently, recentering with respect to SATT (or SATC) can yield a prediction interval (PI) for SATT (or SATC) that is substantially different from a confidence interval (CI) for SATE. It follows that inference on SATE has incorrect coverage and/or is inefficient for the estimation of SATT (or SATC).\footnote{This holds even though SATE, SATT, and SATC are equal in expectation and the difference-in-means is used to estimate all three estimands. Note that under the super-population model, there is no difference between the equivalent three population estimands, as the Population Average Treatment Effect (PATE) and the Population Average Treatment Effect of the Treated (PATT) are equal when the treatment assignment is done at random.} The key result of the paper is the derivation of valid inference on a new and general class of causal estimands that results from any mixing between SATT and SATC. In addition, we present a plug-in estimator for the mixing between these two estimands that maximizes the accuracy -- minimizes the variance -- of the recentered difference-in-means estimator. We call this causal estimand the Sample Average Treatment Effect Optimal (SATO). SATE is a specific case of this general class in which the weighting between SATT and SATC is done by the probability of treatment assignment. SATO is also equal to SATE under a constant/homogeneous treatment effect model.\footnote{A related literature discusses the idea of an optimal estimand in terms of covariate balance in observational studies \citep{Crump_et_al2009, Li2016}. \cite{Crump_et_al2009} suggested a procedure for choosing the optimal estimand in observational studies where there is limited overlap in the covariates. The population overlap issue does not arise in randomized experiments.}  

The objective of the paper is to develop limiting distribution results for inference on sample average treatment effects (e.g., SATE, SATT). \cite{Imbens2004} made explicit the difference between SATE/PATE and SATT/PATT regarding inference, and noted that ``[w]ithout further assumptions, the sample contains no information about the PATE beyond the SATE.'' Therefore, sample estimands are often of direct interest because one does not wish to make additional assumptions.
Inference on SATO yields a PI that has correct coverage (of SATO) and is substantially shorter, more accurate, than a CI for SATE. The heterogeneity in the response of different units to the treatment is the driving force behind the gain in accuracy. The change of estimand yields a recentered difference-in-means test statistic that is sensitive to differences in the variance of the treated and control units, and not only to the mean impact of the treatment. Our approach detects treatment effect heterogeneity by changing the estimand (and corresponding variance formulas), while keeping fixed the test statistic. Detecting and testing for heterogeneity in the treatment effect is of increasing importance to applied researchers \citep{peng2015}.

Furthermore, in some applications the estimand of interest is SATT and using a CI for SATE is inaccurate. An example of when SATT can be the estimand of interest is an attributable treatment effect model in which the treatment effect varies across units \citep{rosenbaum2001, feng2014, keele2017}. In general, PIs for SATT are not guaranteed to have correct coverage of SATE. We provide analytical and simulation based evidence for when a PI for SATT has approximately correct coverage of SATE and when it does not. The key factor behind the differences is the variance of the treatment effect distribution. As the impacts of the treatment are more heterogeneous, inference on SATE differs from inference on SATO.

Our possibly surprising results do have an intuitive interpretation: unlike in the case of the super-population model in which PATE equals PATT and PATC, in the case of Neyman's finite sample model, the three estimands differ from each other and each recentering choice yields a different variance expression for the limiting distribution of the recentered difference-in-means.
Accuracy differences between inference on SATT relative to inference on SATE come from two channels. First, in the case of SATT (and SATC), one does not need to use conservative bounds for the variance estimator as it is point-identified. Second, the change of estimand from SATE to SATT (or SATC) changes the variance of the recentered difference-in-means. We discuss and decompose the conditions for accuracy gains from each one of these channels.  
\cite{robins1988} studied this phenomenon for the special case of binary outcomes. He emphasized that a CI for SATE does not yield correct coverage of SATT (or STAC). We extend his result in various ways, including by providing results for non-binary outcomes and providing conditions under which the SATT and SATC prediction intervals yield gains in accuracy. Unlike Robins, we discuss inference on possible mixes between SATT and SATC and provide a criterion for a plug-in estimator for the optimal mixing between the two estimands.

Our results also extend \cite{rigdon2015}, who showed that PIs for SATT (and SATC) can be combined to construct a CI for SATE in the context of binary outcomes. We generalize their result to the case of non-binary outcomes. In addition, we show that when using the difference-in-means test statistic, combining PIs for SATT and SATC yields the exact same conservative CI for SATE as one would have gotten by using a bound for the variance of the difference-in-means, which assumes that the correlation between potential outcomes is one. More efficient CIs can be constructed directly using sharper bounds on the unobserved correlation between potential outcomes \citep{aronow2014sharp}. Taken together, our results prove that as long as the test statistic is the difference-in-means, combining separate PIs for SATT and SATC, which have been derived in a non-conservative procedure, does \emph{not} allow one to construct a CI for SATE that is more efficient than existing procedures that use conservative variance estimators and directly conduct inference on SATE.

Inference on SATT and SATC may also be of interest even when the primary parameter of interest is SATE. SATT represents the effect of the treatment on the units that have been assigned/exposed to it, and SATC represents the effect of the treatment on the units that have not been exposed to it under the observed treatment allocation. To infer whether the treatment \emph{had} an effect on average, the key parameter of interest is SATT. However, in order to infer whether treatment \emph{will} have an average effect if applied to all the sample units, information on both SATT and SATC is needed. The difference in PIs for each of these estimands is a diagnostic tool for understanding how much of the uncertainty about SATE arises from uncertainty about SATT and how much of it arises from our uncertainty about whether the treatment would have had an effect on the controls (SATC). Decomposing the uncertainty over the SATE can be used as a diagnostic tool for understanding the effects of an intervention, especially in medical trials where heterogeneity in treatment effects is common.   

The remainder of this paper is organized as follows. Section (\ref{sec: real data motivation example}) presents an illustrative data example. Section \eqref{sec: Setting, definitions, and notation} describes the theoretical framework, definitions, and notation that are used throughout the paper. Section \eqref{sec: Theoretical results} describes the key theoretical results of the paper. Section \eqref{sec: Analytical and Monte-Carlo simulations} provides Monte Carlo simulations from three different data generating processes. Section \eqref{sec: Super-population model} extends some of the theoretical results to the super-population sampling model. Section (\ref{sec: Comments and remark}) mentions several notes and remarks about the theoretical results and possible extensions. Section \eqref{sec: Main real the application} discusses the main empirical data application that consists of hundreds of online experiments with millions of observations. Section (\ref{sec: Discussion}) concludes.

\section{Illustrative example} \label{sec: real data motivation example}

Following \cite{rosenbaum2001}, we use data from \cite{tunca1996} who studied the effect of benzene exposure on cytogenetic changes over time. Rosenbaum analyzed this data under an attributable treatment effect model in which the estimand of interest is SATT. 
The treated group ($m=50$) is shoe workers in the area of Bursa, Turkey, and the control group ($N-m = 20$) is residents of Bursa area who were not exposed to benzene. Although the treatment was not assigned at random because this is an observational study, for illustrative purposes we follow previous researchers and assume that the treatment assignment was unconfounded. 

\begin{figure}[ht]
\centering
\caption{ \textbf{Confidence/prediction intervals for average treatment effects \citep{tunca1996}} }
\label{fig: Rosenbaum (2001) CI}
\includegraphics[scale=0.8]{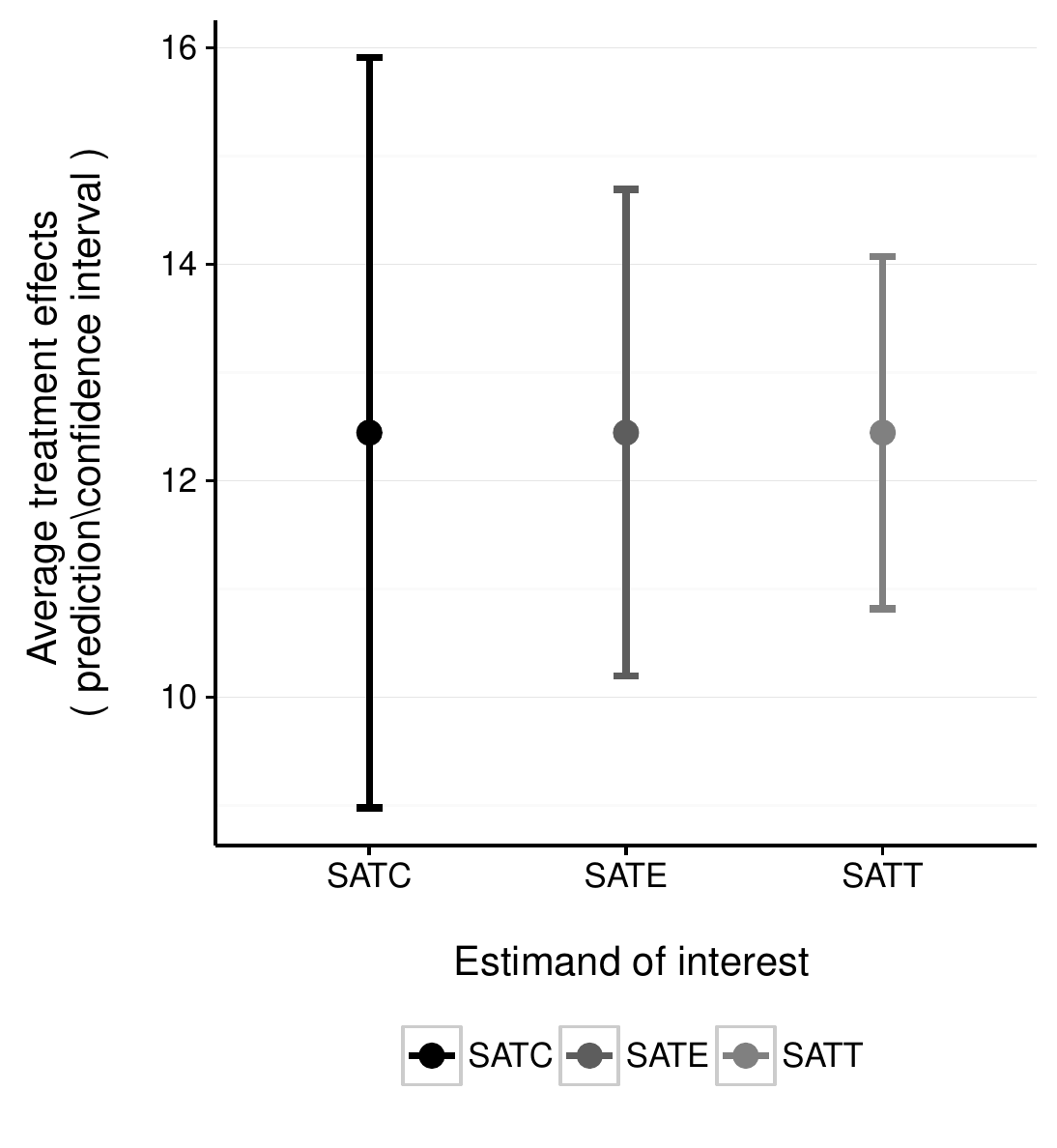}
\end{figure}

Figure \eqref{fig: Rosenbaum (2001) CI} shows a CI (SATE) and PIs (SATT and SATC) for the sample average treatment effects. The PI for SATT is $27\%$ shorter, than a CI for the SATE, while the PI for SATC is $55\%$ longer. The figure demonstrates that in this example SATT is more accurately estimated than SATE, which has a meaningful interpretation. It implies that there is less uncertainty about whether the treatment had an effect, and more uncertainty about the effect the treatment will have if expanded to the units that are currently under the control regime.\footnote{Note that, SATE is a weighted average of SATT and SATC by the probability of treatment assignment.} 

In this example if the researcher is interested in SATT, then a CI for SATE will be overly conservative and inefficient. However, if the estimand of interest is SATC, then a CI for SATE will have incorrect coverage and will be too short; i.e., it will over reject the null hypothesis that SATC is equal to zero.

In what follows we derive analytical results for inference on any weighting between SATT and SATC, which SATE is a special case of. The optimal weighting, in terms of accuracy, between the two estimands depends on the variance of the units under the control and treatment regimes and the correlation between potential outcomes. In the simulations and the main empirical application we show how to conduct inference on an estimand that is a mix between SATT and SATC and that it can detect treatment effect heterogeneity even when SATE is equal to zero.    

\section{Setting, definitions and notation} \label{sec: Setting, definitions, and notation}

We follow Neyman's finite population causal model. Consider a fixed finite population of $N \geq 4$ units and two dimensions, $Y(0)$ and $Y(1)$:
\begin{align}
\Pi_N &=  \left\{ ( Y(0)_{1N}, Y(1)_{1N} ), \; 
( Y(0)_{2N}, Y(1)_{2N} ), \dots, \;
( Y(0)_{NN}, Y(1)_{NN} )  \right\}
\end{align}  
The researcher observes a random sample of $m$ units from the finite population:
\begin{align}
\Pi_N^1 = \left\{ Y(1)_{1N}, Y(1)_{2N}, \dots, Y(1)_{NN}  \right\} 
\end{align} 
and the vector of treatment indicators, $\bm{T} = \left(T_1,\dots,T_N \right)$, represents the $m$ units that are randomly sampled to the treatment group.\footnote{We study the classic Neyman causal model for a finite population. For a review of the framework and the classic CLT results under the finite population model, see \cite{peng2016}. Note that, the randomization model implies that the number of treated units, $m$, is a fixed number and not a random variable. The only random component is the treatment indicators. } The remaining $N-m$ units are assigned to the control group and they form a random sample of $N-m$ units from the finite population $\Pi_N^0 = \left\{ Y(0)_{1N}, Y(0)_{2N}, \dots, Y(0)_{NN}  \right\}$, which is represented by the vector of indicators, $ \left(1-T_1,\dots,1-T_N \right)$.

The probability of each unit being assigned to the treatment regime is $p = \frac{m}{N}$ and the correlation between the treatment assignment of each two units is negative:
\begin{align}
\cov{T_i,T_j} &= p \cdot \left( \frac{m-1}{N-1} - p \right) < 0
\end{align}  
Let $Y_i$ be the outcome of interest for unit $i$, which is a function of $Y(\bm{T})$. We assume SUTVA \citep{holland1986} is satisfied: 
\begin{align}
Y_i(\bm{T})=Y_i(T_i)
\end{align}
and $T_i$ is assigned at random:
\begin{align}
Y(1), Y(0) \independent T
\end{align}

Let $\tau_i$ denote the treatment effect on unit $i$, $Y_i(1) - Y_i(0)$, and let the vector of treatment effects in the finite population be denoted by $\bm{\tau} = \bm{Y(1)} - \bm{Y(0)}$. 
The researcher might be interested in conducting inference on several possible average treatment effect estimands: 
\begin{align}
\text{SATE} &= \frac{1}{N} \cdot \sum_{i=1}^N \tau_i  
\\
\text{SATT} &= \frac{1}{m} \cdot \sum_{i=1}^N \tau_i \cdot T_i, \quad \text{where} \quad m = \sum_{i=1}^N T_i 
\\
\text{SATC} &= \frac{1}{N-m} \cdot \sum_{i=1}^N \tau_i \cdot (1-T_i)
\end{align}     
We do not impose any parametric assumptions on the relationship between $Y(1)$ and $Y(0)$, and allow $\tau$ to vary across units. 
It is important to note that unlike SATE, both SATT and SATC are random variables, even conditional on the sample, as they are a function of $\bm{T}$. 
Because SATT is a random variable, a random interval that contains SATT a with probability of $1-\alpha$ is usually referred to as a prediction interval (PI) or, in Bayesian terminology, a credible interval, rather than a CI. We are not the first ones to discuss inference on a causal estimand that is a random variable \citep{robins1988, hansen2009}. 

\section{Theory} \label{sec: Theoretical results}
\label{sec: Theory fixed/finite population}

In this section we discuss how to conduct inference, without parametric or distributional assumptions, on a general class of average treatment effect estimands:
\begin{align}
\label{eq: general class of estimands}
\omega \att + (1-\omega) \atc
\end{align}

We first present theoretical results on inference about SATT (or SATC) that generalize \cite{robins1988} to non-binary and continuous outcomes. Next we conduct a comparison between inference on SATE and SATT in terms of accuracy and coverage (i.e., type I error). We show that a PI for SATT (or SATC) can be substantially different (smaller or larger) than a CI for SATE, and we decompose the differences to those that result from a change of estimand and those that are the result of being forced to use a conservative variance estimator when conducting inference on SATE. At the end of this section we derive the accuracy-maximizing estimand. This is the mix between SATT and SATC that can be estimated with the highest level of precision while holding fixed the test statistic. SATE is a specific case of the general class of estimands described in equation \eqref{eq: general class of estimands} when $\omega = p$.  

\subsection{Inference on SATT (and SATC) relative to SATE}
The classic estimator for SATE (and SATT) is the difference-in-means between the outcomes of units under the two treatment regimes. For notational convenience, denote the difference-in-means estimator by $\td$:
\begin{align}
\label{eq: difference in means definition}
\td \equiv \frac{1}{m} \sum_{i=1}^N Y_i \cdot T_i - \frac{1}{N-m} \sum_{i=1}^N Y_i \cdot (1-T_i)
\end{align}

Lemma \eqref{lemma: decomposition of the difference in means} shows that the difference-in-means can be decomposed to three terms of which only two depends on $\bm{\tau}$. In equation \eqref{eq:tdiff_satt}, the first expression is a function of $Y_i(0)$ and $T_i$ and is a random variable, the second is a function of $Y_i(0)$ and is not a random quantity, and the third is SATT. 
Lemma \eqref{lemma: decomposition of the difference in means} motivates the use of $\td$ for estimating SATT (or SATC). When recentering $\td$ w.r.t. SATE the variance of the recentered estimator does not change; however, recentering w.r.t. SATT (or SATC) does change the variance calculations. This raises the question of whether inference on SATT (or SATC) has correct coverage of SATE and vice versa. Next we address this question both analytically and using Monte Carlo simulations.
\begin{lemma}
\label{lemma: decomposition of the difference in means}
The difference-in-means test statistic $\td$ can be decomposed into three expressions:
\begin{align}
\label{eq:tdiff_satt}
\diff{\bm{Y},\bm{T}} &= \frac{N}{m \cdot (N-m)} \cdot \sum_{i=1}^N Y_i(0) \cdot T_i -\frac{1}{N-m} \cdot \sum_{i=1}^N Y_i(0) + \text{SATT}
\end{align}
or
\begin{align*}
\diff{\bm{Y},\bm{T}} &= \frac{N}{m \cdot (N-m)} \cdot \sum_{i=1}^N Y_i(1) \cdot T_i -\frac{1}{N-m} \cdot \sum_{i=1}^N Y_i(1) + \text{SATC}
\end{align*}
See Appendix \eqref{proof: diff-in-means decomposition} for the proof. 
\end{lemma}

Lemma \eqref{lemma: variances of recentered diff-in-means} describes the variance of $\td$ when it is recentered with respect to different \emph{causal} estimands of interest. The variance of $\td - \text{SATE}$ contains the parameter $\rho$, which cannot be observed and must be bound when conducting inference. Lemma \eqref{lemma: decomposition of the difference in means} illustrates how by changing the estimand from SATE to SATT (or SATC), still using the difference-in-means test statistic, it is possible to conduct inference without needing to either know or bound $\rho$.  
\begin{lemma}
\label{lemma: variances of recentered diff-in-means}
The variance of the difference-in-means when recentered w.r.t. SATE, SATT or SATC is 
\begin{align*}
\var{ \td - \text{SATE} } &= 
\frac{1}{N \cdot (1-p) \cdot p} \cdot \left[ p^2 \cdot \sigma_0^2 + (1-p)^2 \cdot \sigma^2_1
+ 2 p (1-p) \cdot \rho \cdot \sigma_0 \cdot  \sigma_1 \right]
\\
\var{ \td - \text{SATT} } &= \frac{1}{p \cdot (1-p) \cdot N} \cdot \sigma_0^2 
\\
\var{ \td - \text{SATC} } &= \frac{1}{p \cdot (1-p) \cdot N} \cdot \sigma_1^2
\end{align*}
where $\rho \equiv \frac{\cov{ Y(1), Y(0) }}{ \sigma_0 \cdot \sigma_1 } $ is the correlation between the potential outcomes under the treatment and control regimes, and $\sigma^2_j$ is defined as $\sigma^2_j = \frac{ \sum_{i=1}^N \left( Y(j)_i - \bar{Y}(j) \right)^2 }{N-1}$.
See Appendix \eqref{proof: variance of different recentered diff-in-means} for the proof. 
\end{lemma}

The $\td - \text{SATT}$ can be more accurately estimated relative to the $\td - \text{SATE}$ when $\frac{\sigma_1}{\sigma_0} > 1$, $\rho$ is sufficiently high, and $p$ is not too high (e.g., $p=1/2$). Theorem \eqref{theorem: gains from estimand} shows that regardless of the value of $\frac{\sigma_1}{\sigma_0}$ there is a threshold level of $\rho$ that below it $\var{ \td - \text{SATE} } \leq \var{ \td - \text{SATT} }$ and above which $\var{ \td - \text{SATE} } > \var{ \td - \text{SATT} }$. Notice that according to Theorem \eqref{theorem: gains from estimand}, it is simple to empirically test whether $\bar{\rho}$ is negative. All that is needed is to conduct a one-sided hypothesis test of the null: 
\[
H_0: \; \frac{\sigma_1}{\sigma_0} \leq \sqrt{\frac{1-p^2}{(1-p)^2}}
\] 
and if the null hypothesis is rejected we can infer that $\bar{\rho}<0$.  
\begin{theorem}
\label{theorem: gains from estimand}
For all $\sigma_0$ and $\sigma_1$ such that $\sigma_0 < \sigma_1$:
\begin{enumerate}
\item  There exists a threshold level of $\rho$, $\bar{\rho}$ such that
\begin{align*}
\rho \leq  \bar{\rho} \Rightarrow \var{ \td - \text{SATE} } \leq \var{ \td - \text{SATT} }
\\
\rho >  \bar{\rho} \Rightarrow \var{ \td - \text{SATE} } > \var{ \td - \text{SATT} }
\end{align*}
\item When $\frac{\sigma_1}{\sigma_0} > \sqrt{\frac{1-p^2}{(1-p)^2}}$ then, $\bar{\rho} < 0$.
\end{enumerate}
See Appendix \eqref{proof: gains from estimand change} for the proof.
\end{theorem}

In practice, the correlation between potential outcomes is not observed and to estimate $\var{ \td - \text{SATE} }$ it is required to use a bound for $\rho$. The most commonly used estimator was proposed by Neyman and it ignores the correlation component all together. It can be re-written as:
\begin{align}
\mathbb{V}_{\text{Neyman}} &= 
\frac{1}{m} \sigma_1^2 + \frac{1}{N-m} \sigma_0^2
%
\nonumber \\ &=
\frac{1}{N \cdot p \cdot (1-p)} \left( \sigma_1^2 (1-p) + \sigma_0^2 p \right).
\end{align}
This variance estimator is a consistent estimator for the variance of the difference-in-means under the super-population sampling model, and it can be used to conduct inference on the Population Average Treatment Effect (PATE). It also corresponds to Neyman's second variance estimator when the estimand is SATE \citep{neyman1923}. A less conservative variance estimator for the difference-in-means bounds $\rho$ at 1:
\begin{align}
\label{eq: corr=1 variance estimator}
\mathbb{V}_{\rho=1} &= 
\frac{1}{m} \sigma_1^2 + \frac{1}{N-m} \sigma_0^2 - \frac{(\sigma_1 - \sigma_0)^2}{n}
\nonumber \\ &=
\frac{1}{N \cdot (1-p) \cdot p} \left( p \cdot \sigma_0 + (1-p) \cdot \sigma_1 \right)^2
\end{align}

\begin{theorem}{(Asymptotic distribution of $\td-\att$ )}
\label{theorem: CLT - fixed sample}
The difference-in-means recentered w.r.t. SATT follows a standard Normal distribution under two regularity conditions. When the following two conditions are satisfied:
\begin{align}
N-m \rightarrow \infty, \quad m \rightarrow \infty
\quad
\text{and}
\quad
\frac{ \underset{ 1 \leq i \leq N }{\max}  \left( Y(0)_{Ni} - \bar{Y}(0)_{N}  \right)^2 }{ \sum_{i=1}^N \left( Y(0)_{Ni} - \bar{Y}(0)_{N}  \right)^2  }
\cdot
\max \left( \frac{N-m}{m}, \frac{m}{N-m}  \right)
\rightarrow 0 
\nonumber
\end{align}
then
\begin{align}
\frac{ \td - \text{SATT}  }{ \sqrt{ \var{ \td - \text{SATT} } } } \overset{d}{\rightarrow} N(0,1)
\end{align}
see Appendix \eqref{proof: CLT - fixed sample} for the proof.\footnote{These are not the weakest possible conditions.} 
\end{theorem}
\noindent Therefore, a $1-\alpha$ prediction interval for SATT is
\begin{align}
\left[ \td - z_{1-\alpha/2} \cdot \hat{\sigma}_0 \cdot \sqrt{k(N,m)}, \;  \td + z_{1-\alpha/2} \cdot \hat{\sigma}_0 \cdot \sqrt{k(N,m)}   \right]
\end{align}
where
\begin{align*}
k(N,m) &= \frac{1}{p \cdot (1-p) \cdot N}
\end{align*}

Theorem \eqref{theorem: CLT - fixed sample} establishes that the limiting distribution of the $\td - \text{SATT}$, when standardized using the variance formulas in Lemma \eqref{lemma: variances of recentered diff-in-means}, is standard Normal. An equivalent derivation of the limiting distribution can be carried out for SATC. The theorem for $\td - \text{SATC}$ follows immediately using an analog proof. 

\cite{rigdon2015} showed how to derive a CI for SATE by combining two PIs for SATT and SATC. They focused on binary outcomes for which PIs have been derived using past theoretical results \citep{robins1988, rosenbaum2001}. Theorem \eqref{theorem: Combining PIs for CI for SATE - corr=1} provides two key results. First, it is a generalization of Theorem (1) in \cite{rigdon2015} and shows how a CI for SATE can be constructed in any randomized control trial regardless of whether the outcomes are binary or continuous. This part of the theorem follows directly from the derivations in \cite{rigdon2015}. Second, it shows that combining PIs for SATT and SATC using a Bonferroni-type adjustment, as was done by Rigdon and Hudgens, yields \emph{exactly} the same CI as when one constructs a CI for SATE directly and uses a conservative variance estimator that bounds $rho$ at 1, $\mathbb{V}_{\rho=1}$. Theorem \eqref{theorem: Combining PIs for CI for SATE - corr=1} implies that combining PIs for SATT and SATC yields a more conservative CI than using a CI based on a sharper bound for $\rho$, such as that of \cite{aronow2014sharp}.     

\begin{theorem}
\label{theorem: Combining PIs for CI for SATE - corr=1}
Let $[L_{\att}, U_{\att}]$ and $[L_{\atc}, U_{\atc}]$ be PIs for SATT and SATC according to the variance formulas in Lemma \eqref{lemma: variances of recentered diff-in-means}; then a CI for SATE is
\[
\left[  p \cdot L_{\att}  + (1-p) \cdot L_{\atc}, \quad p \cdot U_{\att}  + (1-p) \cdot U_{\atc}  \right]
\]
and is equal to
\[
\left[ \td - z_{1-\alpha/2} \cdot \sqrt{\hat{\mathbb{V}}_{\rho=1}}, \;  \td + z_{1-\alpha/2} \cdot \sqrt{\hat{\mathbb{V}}_{\rho=1}}  \right]
\]
See Appendix \eqref{proof: Combining PIs for CI for SATE - corr=1} for the proof.   
\end{theorem}

\begin{theorem}
\label{theorem: prediction vs. confidence intervals - fixed sample}
When $\sigma_1 \neq \sigma_0$, a PI for either SATT or SATC is shorter than a CI for SATE that uses either $\mathbb{V}_{\text{Neyman}}$ or $\mathbb{V}_{\rho=1}$. \\
See Appendix \eqref{proof: prediction vs. confidence intervals - fixed sample} for the proof.   
\end{theorem}
According to Theorem \eqref{theorem: prediction vs. confidence intervals - fixed sample}, whenever $\sigma_0 < \sigma_1$ ($\sigma_0 > \sigma_1$), a PI for SATT (SATC) is shorter than a CI for SATE using Neyman's variance estimator. The gain in terms of interval length (in \%) is
\begin{align}
1 -\frac{ 1 }{ \sqrt{ \left( \frac{\sigma_1^2}{\sigma_0^2} (1-p) + p \right) }  } 
\end{align}
and the gains are decreasing with respect to $p$.

\begin{figure}[ht]
\centering
\caption{ \textbf{Percentage gains from a PI for SATT instead of a CI for SATE} }
\label{fig: gain analytical SATT vs SATE}
\includegraphics[scale=0.8]{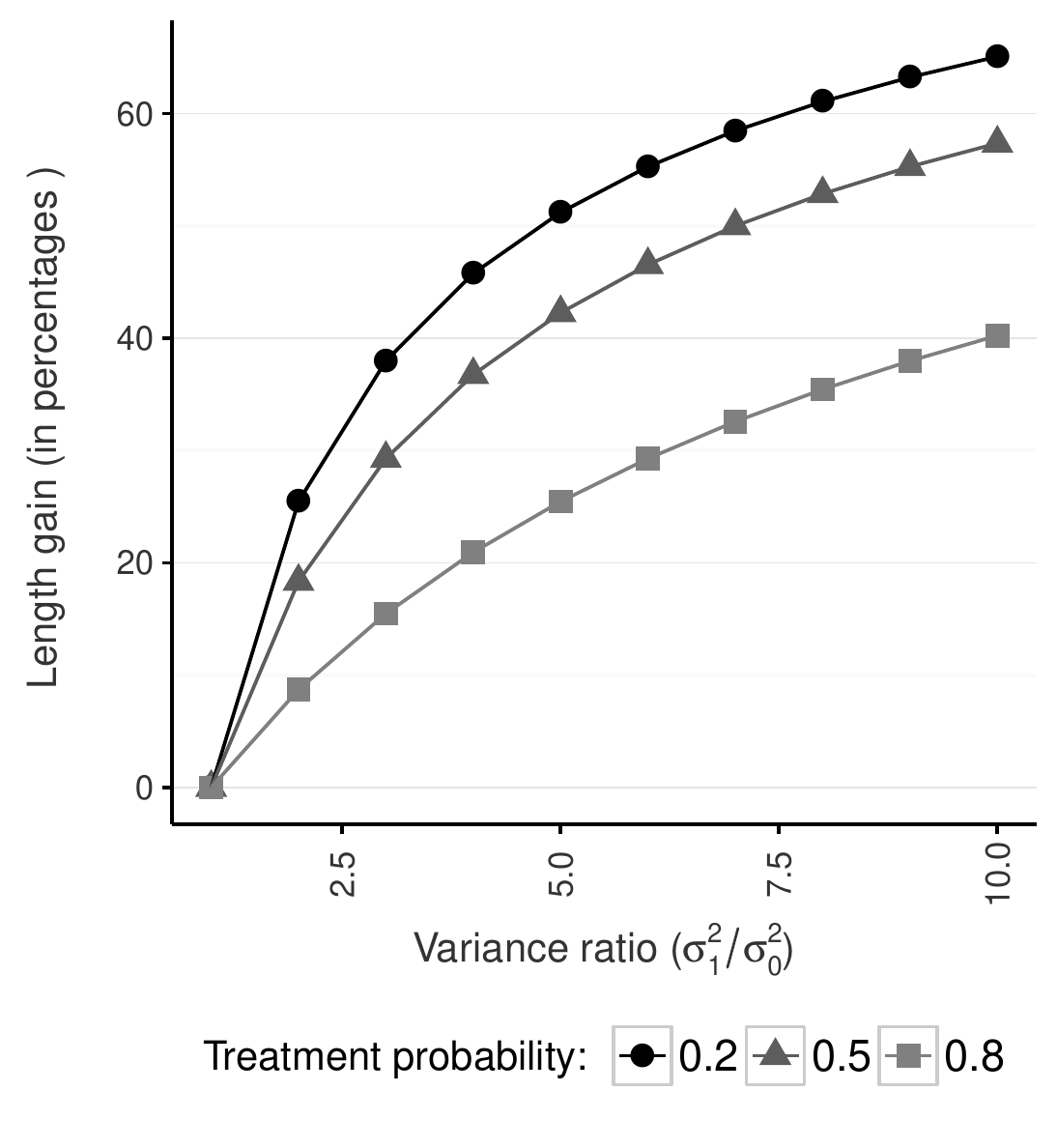}
\end{figure}

The interval length gains (in percentages) of a PI for SATT relative to a CI for SATE are illustrated in Figure \eqref{fig: gain analytical SATT vs SATE}. In a balanced design, $p = \frac{1}{2}$, with a variance ratio of two, the length/accuracy gain is 18.35\% relative to a CI for SATE that is based on Neyman's variance estimator. Figure \eqref{fig: gain analytical SATT vs SATE} shows the accuracy gains from estimating SATT instead of SATE for different levels of $p$, which can be substantial.  

Note that Theorem \eqref{theorem: prediction vs. confidence intervals - fixed sample} does not cover other estimators for the variance of $\var{ \td - \mathrm{SATE} }$, such as that of \cite{aronow2014sharp}. CIs based on these other variance estimators may be shorter than our PIs even when $\sigma_1 \neq \sigma_0$.  
To address these variance estimators, we derive the accuracy gains for a more general case in which $\rho$ can be bounded by $\rho^*$, $\rho \leq \rho^*$. As the variance of $\var{ \td - \mathrm{SATE} }$ is increasing w.r.t. $\rho$, it follows that substituting $\rho^*$ with $\rho$ yields a conservative variance estimator that is smaller than Neyman's variance estimator. The idea of substituting a bound of $\rho$ instead of the true parameter value was proposed before in the literature \citep{reichardt_1999, aronow2014sharp}.  
The percentage gain in terms of CI length is
\begin{align*}
1 - \frac{1}{ \sqrt{ p^2 + (1-p)^2 \cdot \left( \frac{\sigma_1}{\sigma_0} \right)^2 + 2p(1-p) \cdot \rho^* \cdot \cdot \frac{\sigma_1}{\sigma_0} } } 
\end{align*}

To obtain intuition about how recentering the difference-in-means w.r.t. SATT provides an accuracy gain, it is useful to decompose the variance of the difference-in-means. The decomposition in equation \eqref{eq: variance difference-in-means decomposition} shows that when recentering w.r.t. SATT we cancel two of the elements in the variance of the difference-in-means, $\var{ \text{SATT} }$ and 
 $\cov{ \text{SATT}, \; \frac{N}{m (N-m)} \cdot \sum_{i=1}^N Y(0) T_i  }$, and in doing so reduce the uncertainty.    
\scriptsize
\begin{align}
\label{eq: variance difference-in-means decomposition}
\var{ \td - \text{SATE} } &= \var{ \text{SATT} } 
+ \var{\frac{N}{m (N-m)} \cdot \sum_{i=1}^N Y(0) T_i }
+ 2 \cdot \cov{ \text{SATT}, \; \frac{N}{m (N-m)} \cdot \sum_{i=1}^N Y(0) T_i  }
\end{align}
\normalsize

To understand why inference on SATE can provide incorrect coverage of SATT (or SATC), it is useful to look at the mean square distance between the two estimands:
%
\begin{align}
\text{MSE}(\text{SATE}, \; \text{SATT})  &= \frac{1-p}{m} \cdot \sigma^2_{ \tau }
\end{align} 
see Appendix \eqref{proof: derivation of MSE(SATE,SATT)} for the proof. The differences between conducting inference on SATE or SATT does increase with the heterogeneity of the treatment effect. Given a fixed and bounded value of $\sigma_{\tau}^2 < \infty$, those difference will converges to zero as $m \rightarrow \infty$. 

To illustrate how SATT and SATE can differ substantially, consider the following extreme example. There is a finite population of four units (see Table \eqref{tab: example SATE vs SATT}) of which two are assigned to the treatment regime and the other two to the control regime. SATE is exactly zero in this example; however, SATT can be negative, positive, or zero depending on which units are allocated to the treatment group. In this example using inference on SATE (or SATT) to approximate inference on the SATT (or SATE) is misleading. However, inference on SATT can be valuable; and an estimand of interest in its own. 
\begin{table}[ht]
\centering
\caption{ \textbf{Example of when SATT can substantially differ from SATE} }
\label{tab: example SATE vs SATT}
\begin{tabular}{lcc}
\hline \hline 
Unit & $Y(1)$ & $Y(0)$ \\ [0.5ex]
\hline 
1 & 1 & 0 \\
2 & -1 & 0 \\
3 & -100 & 0 \\
4 & 100 & 0 \\
\hline
\end{tabular}
\end{table}

\subsection{Inference a general class of average treatment effects }

In this section we derive inference for the Sample Average Treatment Effect Optimal (SATO), which is the weighting of SATT and SATC that minimizes MSE. The inference is conditional on the sample and uses the difference-in-means test statistics. Formally, SATO is defined as
\begin{align*}
& \ato \equiv \omega^* \cdot \att + (1-\omega^*) \cdot \atc \\
& s.t. \\
& \omega^* = \underset{\omega}{\text{argmin}} \; \var{ \hat{\bar{Y}}_1 - \hat{\bar{Y}}_0 - \ato }
= \underset{\omega}{\text{argmin}} \; \text{MSE} \left( \hat{\bar{Y}}_1 - \hat{\bar{Y}}_0 - \ato,\; 0 \right)
\end{align*}

\begin{lemma}
\label{lemma: SATO variance and omega}
The variance of the difference-in-means when recentered w.r.t. any choice of weighting between SATT and SATC, $\var{ \td - \left( \omega \att + (1-\omega) \atc \right) }$, is
\begin{align*}
 & \frac{1}{Np(1-p)} \left[  p^2 \sigma_0^2 + (1-p)^2 \sigma_1^2 + 2 p(1-p) \cdot \rho \sigma_0 \sigma_1 \right]
\\ & \quad 
+ \frac{\sigma_{\tau}^2}{N-1} \cdot \left[ \omega^2 \cdot \frac{1-p}{p} + (1-\omega)^2 \cdot \frac{p}{1-p} -2 \omega (1-\omega)  \right]
\\ & \quad
- 2 \cdot \left[ \frac{\omega}{(N-1)p} \left[ \sigma_1^2-\rho \sigma_1 \sigma_0   \right] - \frac{1}{N-1} \cdot \sigma_{\tau}^2 + \frac{(1-\omega)}{ (N-1)(1-p) } \left[ \sigma_0^2 - \rho \sigma_1 \sigma_0  \right]  \right]
\end{align*}
See Appendix \eqref{proof: diff-in-means recentered w.r.t SATO variance} for the proof. 
\end{lemma}

Lemma \eqref{lemma: SATO variance and omega} presents the variance of the difference-in-means when it is recentered w.r.t. $\omega \att + (1-\omega)\atc$ for any choice of $\omega$. It follows immediately from Lemma \eqref{lemma: SATO variance and omega}, that the optimal accuracy-maximizing value of $\omega$ is    
\begin{align*}
\omega^* &=  \frac{ \left( \frac{\sigma_1}{\sigma_0} \right)^2 - \rho \cdot \frac{\sigma_1}{\sigma_0}  }{ \left( \frac{\sigma_1}{\sigma_0} \right)^2 + 1 - 2 \rho \left( \frac{\sigma_1}{\sigma_0} \right)}
\end{align*} 

\begin{figure}[ht]
\centering
\caption{ \textbf{Optimal $\omega$ weight of SATT in SATO for different $\frac{\sigma_1}{\sigma_0}$ and $\rho$}}
\label{fig: SATO optimal omega for different variance ratios}
\includegraphics[scale=0.8]{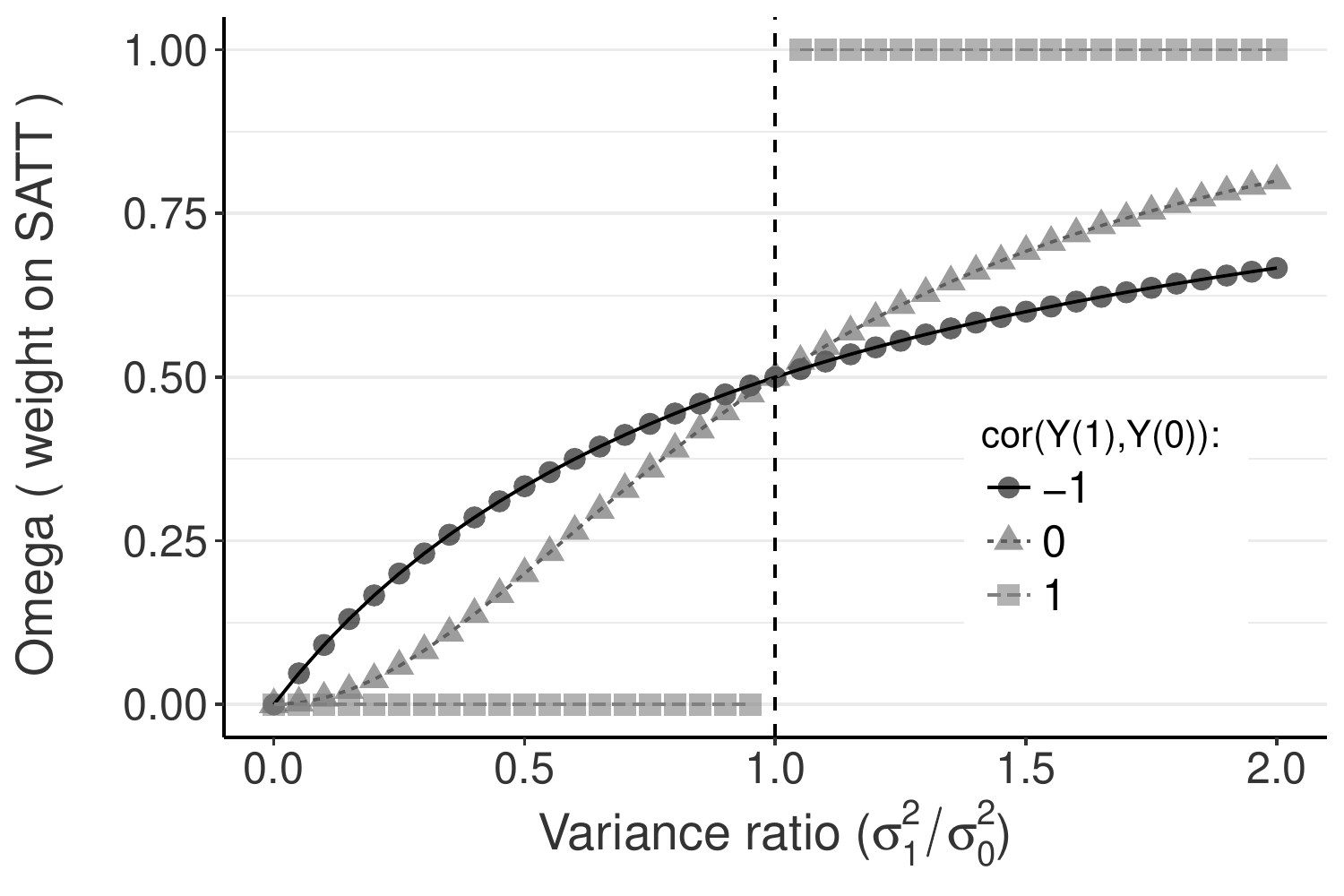}
\begin{minipage}{16cm}
\footnotesize
\emph{Notes}: The above figure describes the relationship between the variance ratio of treated and control units ($\frac{\sigma_1}{\sigma_0}$) and the weight which is assigned to SATT in SATO for three different values of $\rho$. 
\end{minipage}
\end{figure}

The optimal choice of $\omega$ depends on two parameters, $\rho$ and $\frac{\sigma_1}{\sigma_0}$, and is \emph{independent} of $p$. This is in stark contrast to SATE, which is the weighting of SATT and SATC according \emph{only} to the probability of treatment assignment. Figure \eqref{fig: SATO optimal omega for different variance ratios} illustrates how the value of $\omega$ changes depending on $\rho$ and $\frac{\sigma_1}{\sigma_0}$. Two clear patterns stand out. First, the weight that is assigned to SATT is monotonically increasing w.r.t. the variance ratio, $\frac{\sigma_1}{\sigma_0}$. Second, the relationship between $\omega^*$ and $\rho$ is ambiguous and depends on $\frac{\sigma_1}{\sigma_0}$. When $\frac{\sigma_1}{\sigma_0}>0$, as $\uparrow \rho$, the weight assigned to SATT increases, and the opposite is true when $\frac{\sigma_1}{\sigma_0}<0$. 

\begin{figure}[ht]
\centering
\caption{\textbf{The accuracy gain (in percentages) from estimating SATO instead of SATE when $\frac{\sigma_1}{\sigma_0}=1$ for different values of the treatment assignment probability ($p$) and the correlation between potential outcomes ($\rho$)}}
\label{fig: Equal variance SATO gains }
\includegraphics[scale=0.66]{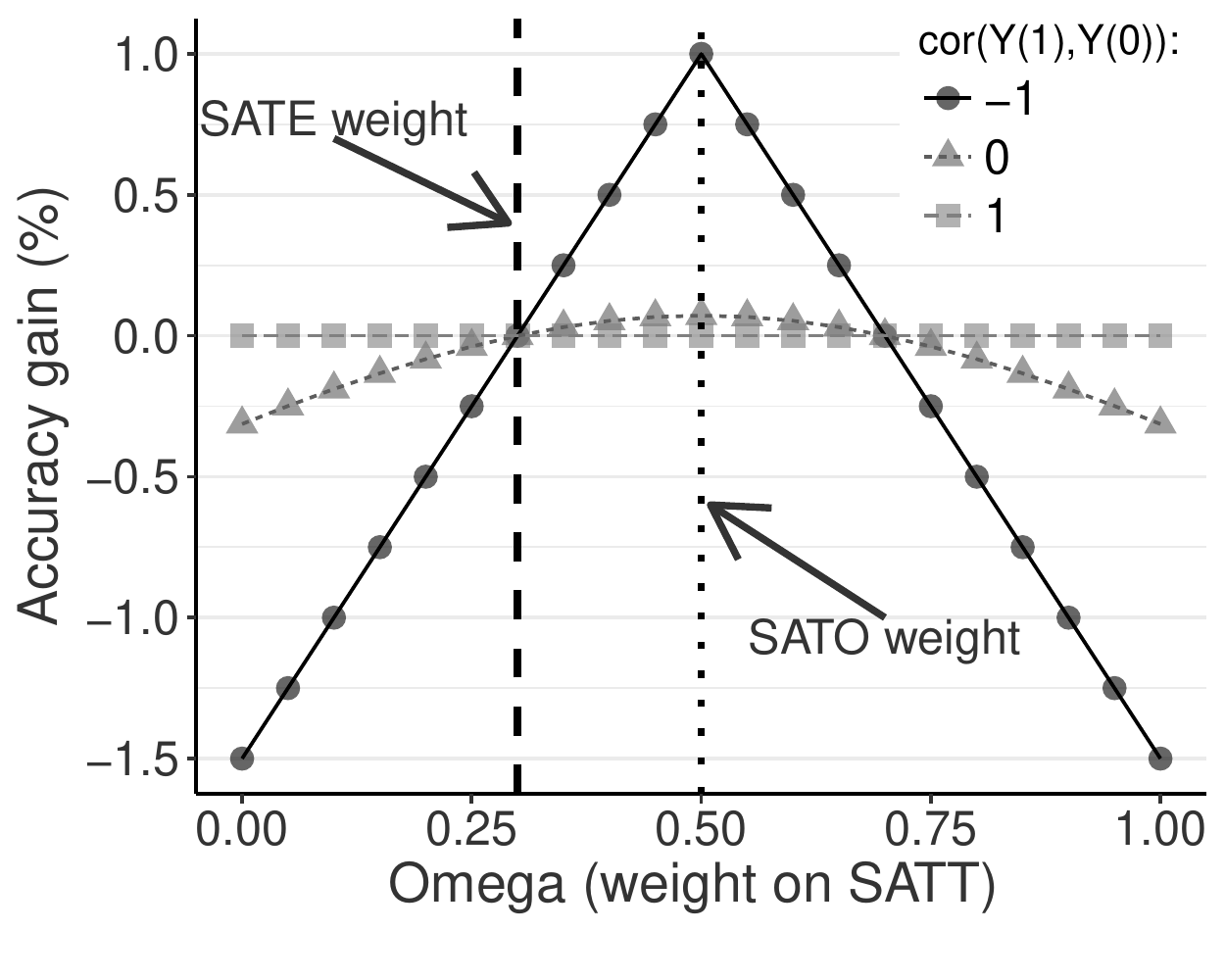}
\includegraphics[scale=0.66]{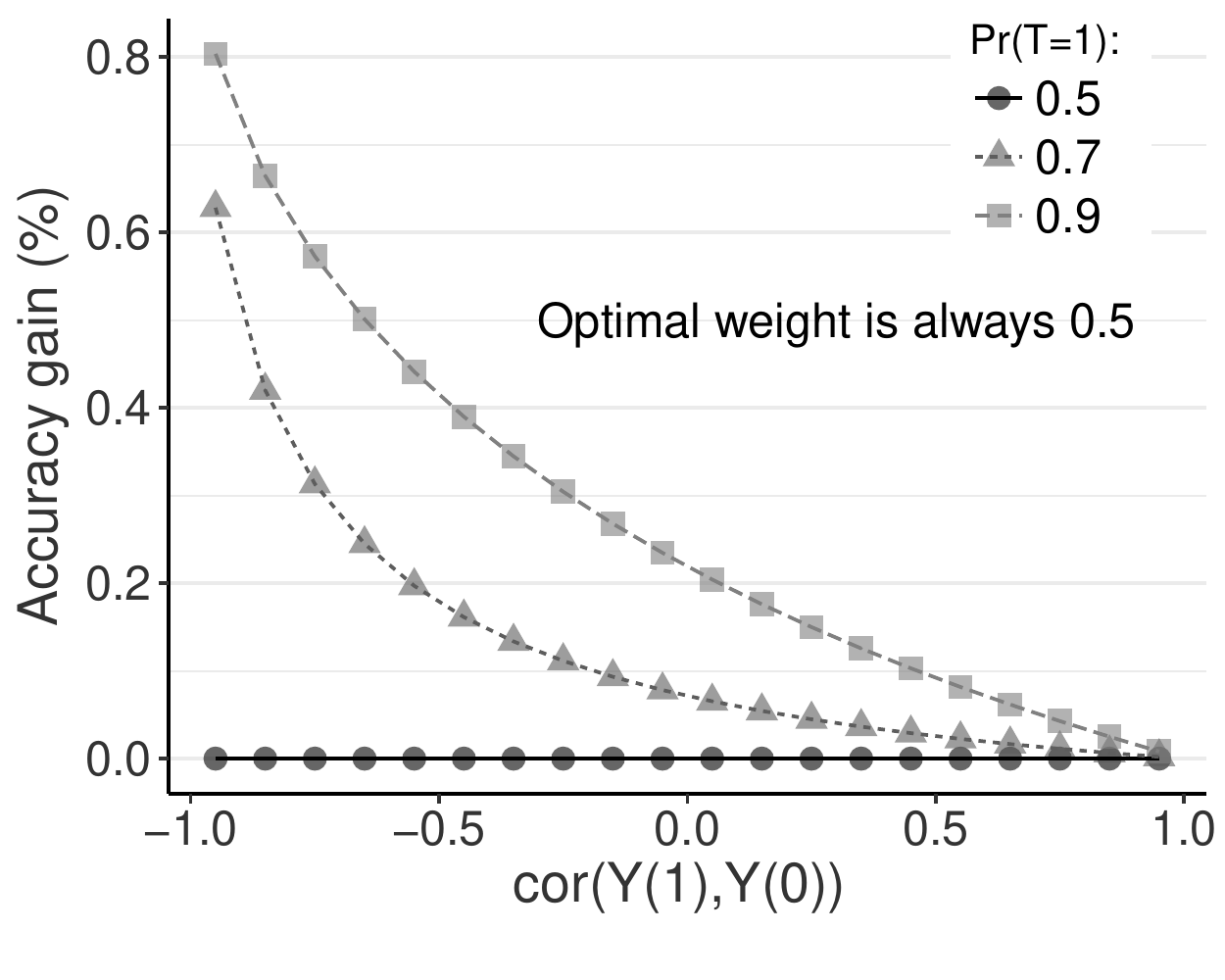}
\begin{minipage}{16cm}
\footnotesize
\emph{Notes}: The left plot describes the relationship between the accuracy gain of conducting inference on a mixing between SATT and SATC relative to SATE as a function of the weight which is assigned to SATT ($\omega$) for three different values of $\rho$. The right plot illustrates the relationship between $\rho$ and the accuracy gains of a PI for SATO relative to a CI for SATE for three different levels of the probability of treatment assignment. In both plots the relationships that are described are for the case in which $\frac{\sigma_1}{\sigma_0} = 0$, which implies that accuracy gains are not the result of a difference in the variance of the outcome among treated units relative to the control units.       
\end{minipage}
\end{figure}

A key question is under what conditions SATE is equal to SATO, which is equivalent to asking whether SATE can ever be the estimand that is estimated most accurately. Two scenarios in which SATO coincides with SATE are (i) a constant treatment effect model and (ii) when $\sigma_1=\sigma_0$ and $p=1/2$.  
When the variance ratio is equal to 1, the optimal weighting between SATT and SATC is half and half ($\omega = 0.5$). Equality of variance between treatment and control units does \emph{not} imply that SATO collapses to SATE. The right plot of Figure \eqref{fig: Equal variance SATO gains } shows the accuracy gain of SATO relative to SATE for different values of $\rho$ and $p$. SATO can be estimated at a higher degree of accuracy for most levels of $\rho$ and for any $p \neq 1/2$, even when $\sigma_1=\sigma_0$. The left plot of Figure \eqref{fig: Equal variance SATO gains } compares the accuracy of inference on SATE to a mixing between SATT and SATC -- not SATO as $\omega^*$ is not used -- for a variety of different $\omega$ weights and three values of $\rho$. SATE can be estimated more accurately than many mixes of SATT and SATC, but it will always -- for a known $\rho$ -- be more noisily estimated than SATO. To conclude, it is hard to justify a weighting of SATT and SATC that will yield the SATE when the objective is to maximize accuracy. SATE is the only mix of SATT and SATC that yields a parameter and \emph{not} a random variable.    

\section{Analytical examples and simulations} \label{sec: Analytical and Monte-Carlo simulations}

\subsection{Additive random coefficients model}
\label{sec: Analytical example - random coefficients}

\begin{figure}[ht]
\centering
\caption{ \textbf{Additive and heterogeneous treatment effect (random coefficient) simulation results: Length and coverage differences} }
\label{fig: random_coefficient simulation}
\includegraphics[scale=0.55]{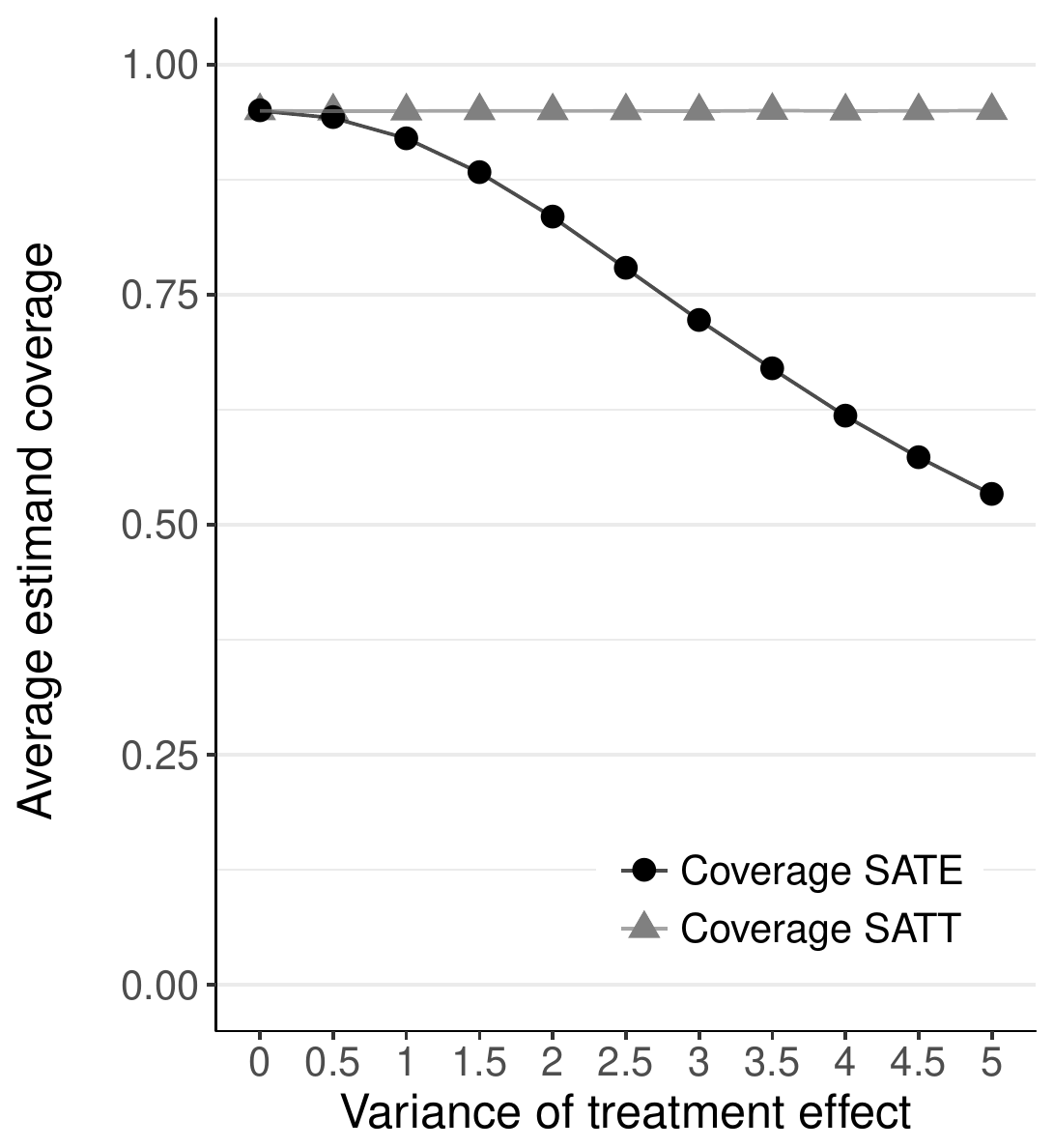}
\includegraphics[scale=0.55]{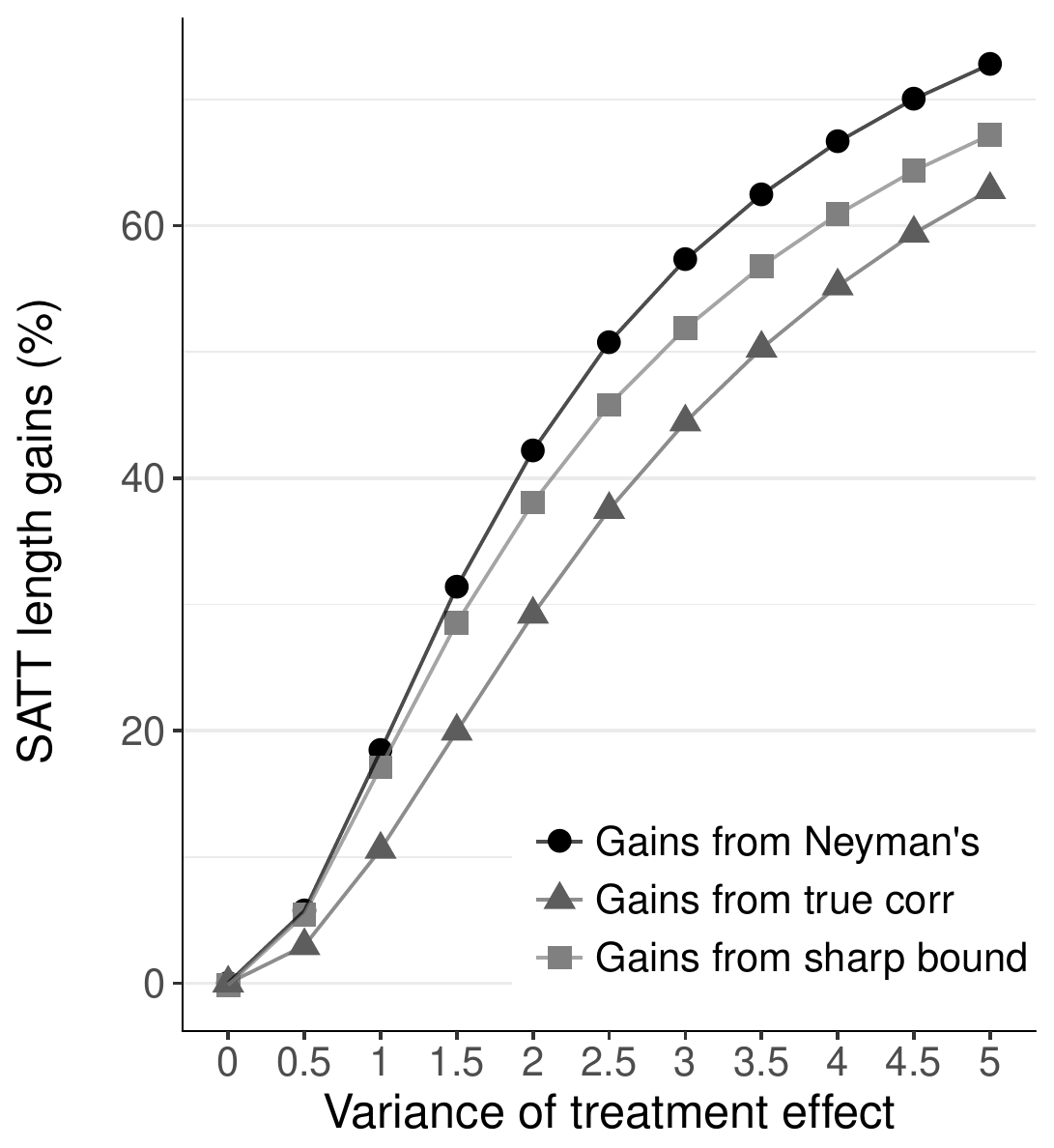}
\begin{minipage}{16cm}
\footnotesize
\emph{Notes}: The left plot shows the coverage of a PI for SATT w.r.t. SATT and SATE. The coverage of SATT has correct size and there is no evidence of over rejection of the null hypothesis. On the other hand, the coverage of SATE becomes worse as the variance of the treatment effect increases. The right plot illustrates the accuracy gains of using inference on SATT relative to SATE. The PI for SATT is compared to CIs for SATE based on three different values of $\rho$. The circles represent the accuracy gains (in percentages) relative to Neyman's variance estimator which assumes that $\rho = 0$. The squares dots use a sharp bound for $\rho$ that have been derived by \cite{aronow2014sharp}. And the triangles compare the length of a PI for SATT to a CI for SATE when the true value of $\rho$ is known to the researcher.   
\end{minipage}
\end{figure}

We consider a simple additive treatment effect model with heterogeneous impacts across units:
\begin{align}
\label{eq: random coefficient model}
Y_i(1) &= \tau_i + Y_i(0)
\end{align}
where $\tau_i$ is a random variable, which for simplicity is assumed not to be correlated with $Y_i(0)$, nor to be correlated with $T_i$ by construction due to the randomization of treatment assignment. The variance of the treated units is larger by $\sigma_{\tau}^2$, and this generates a potential difference between a CI for SATE and a PI for SATT. The variance ratio is:
\[
\frac{\sigma_{\tau}^2}{\sigma_0^2} + 1
\] 
and is increasing with the heterogeneity of the treatment effect. Consider the data-generating process in \eqref{eq: random coefficient dgp} with a sample of $n=1,000$ units. We performed $1,000$ draws of samples and for each one simulated $1,000$ different allocations of the treatment according to $p = 1/2$. 
Figure \eqref{fig: random_coefficient simulation} reports the simulation results. The PIs for SATT are substantially shorter than the CIs for the SATE; this holds both when using Neyman's variance estimator and when using a variance estimator based on a sharp bound for $\rho$, such as that of \cite{aronow2014sharp}.  
The comparison to the CI based on the true $\rho$ parameter shows that the accuracy differences are mainly due to the change of estimand, rather than a conservative variance estimator. This is contrary to the binary outcome example (see below). The simulation is supported by our theoretical results that a change of estimand has larger accuracy impacts as the variance of the treatment effect is higher. The left plot in Figure \eqref{fig: random_coefficient simulation} shows that coverage of the PI for SATT w.r.t. both SATT and SATE. As the treatment effect becomes more heterogeneous the PI for SATT provides worse coverage of SATE and, vice versa, a CI for SATE provides worse coverage of SATT. In a heterogeneous treatment effect environment it is required to use variance formulas that guarantee correct coverage of the estimand of interest, as SATE, SATT, and SATC can differ substantially from one another.    
\begin{align}
\label{eq: random coefficient dgp}
Y_i(0) & \sim N(\mu = 10, \; \sigma_0^2 = 1) \nonumber \\
\tau_i & \sim N(\mu = 0, \; \sigma_\tau^2) \nonumber \\
Y_i(1) &= \tau_i + Y_i(1)
\end{align}

\subsection{Binary outcomes}
\label{sec: Analytical example - binary outcome}

Consider a finite population of size $N$:
\begin{align}
\Pi_N &=  \left\{ ( Y(0)_{1N}, Y(1)_{1N} ), \; 
( Y(0)_{2N}, Y(1)_{2N} ), \dots, \;
( Y(0)_{NN}, Y(1)_{NN} )  \right\}
\end{align}   
where $\left(  Y(0)_{iN}, Y(1)_{N}  \right) \in \lbrace 0, 1 \rbrace^2$, and 
\begin{align*}
p_0 &= \frac{1}{N} \cdot \sum_{i=1}^N Y_{i}(0) 
\quad \text{and} \quad 
p_1 = \frac{1}{N} \cdot \sum_{i=1}^N Y_{i}(1), 
\end{align*}   
The outcome variable is binary and we assume the treatment has a positive effect on average:
\begin{align}
\label{eq: condition analytical example}
& 0 \leq p_0 \leq 1/2 \quad \text{and} \quad p_0<p_1 \leq 1/2
\end{align}
These conditions insure that the positive treatment effect increases the variance of the outcomes among the treated units. 
SATE is equal to $p_1 - p_0$, which implies that when the conditions in equation \eqref{eq: condition analytical example} are met, the variance of the treated units is strictly higher than that of the controls, $p_0(1-p_0) < p_1(1-p_1)$, and the variance ratio is a function of the treatment effect,
\begin{align}
\frac{\sigma_1^2}{\sigma_0^2} &= \frac{(p_0+SATE)(1-p_0-SATE)}{p_0(1-p_0)}
\end{align}
The accuracy gains are increasing with the treatment effect as long as the treatment generates a variance increases among the treated units. Figure \eqref{fig:   binary outcome gain analytical SATT vs SATE} shows the interval length gains (in percentages) from using a PI for SATT in the case of a binary outcome. When $p_0=0.1$ and $p_1=0.2$, the interval for SATT is lower by $23.3\%$ relative to the standard CI that is based on Neyman's variance estimator. 

\begin{figure}[ht]
\centering
\caption{ \textbf{Binary outcome simulation: Length and coverage differences.} }
\label{fig: binary outcome gain analytical SATT vs SATE}
\includegraphics[scale=0.55]{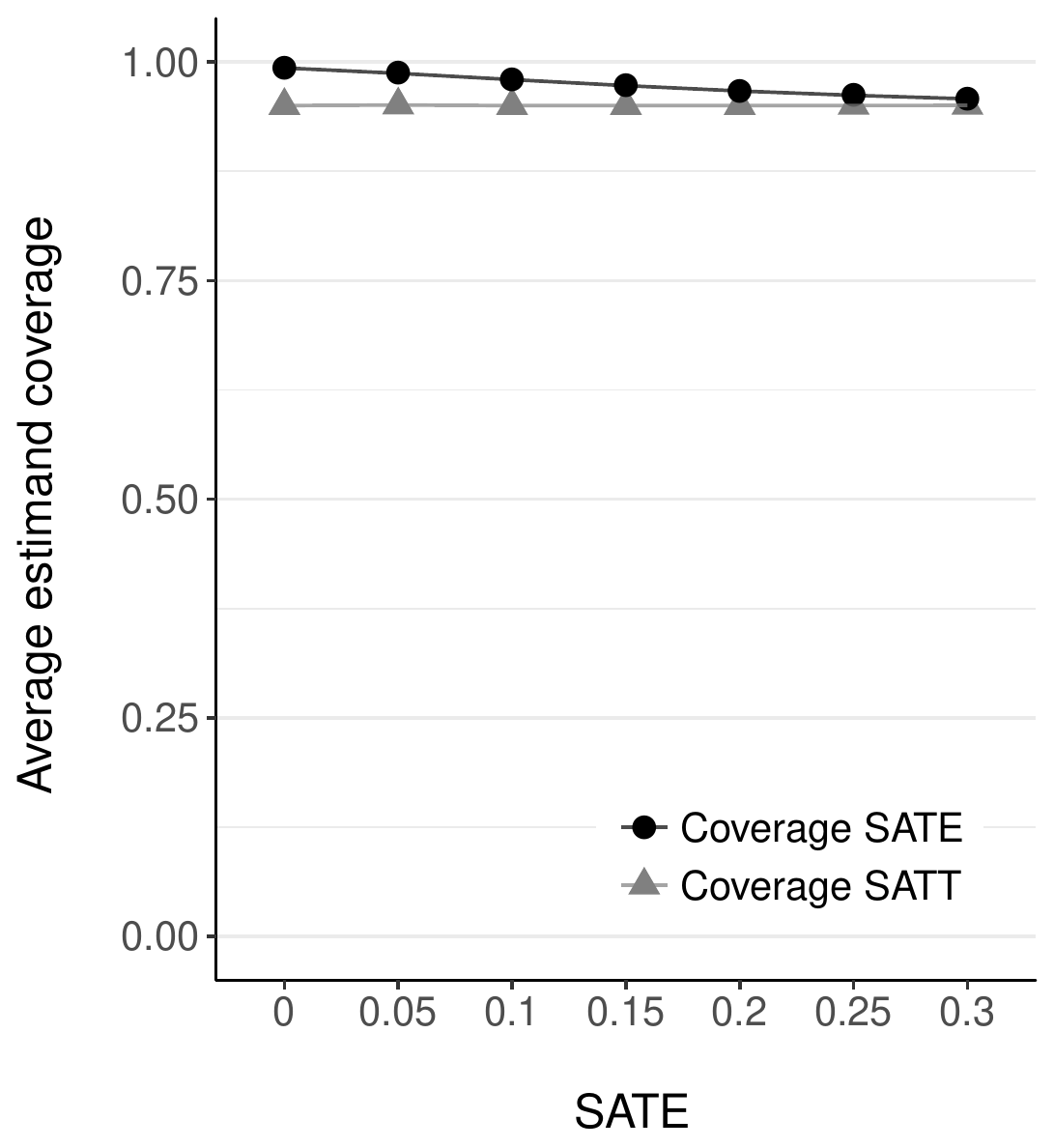}
\includegraphics[scale=0.55]{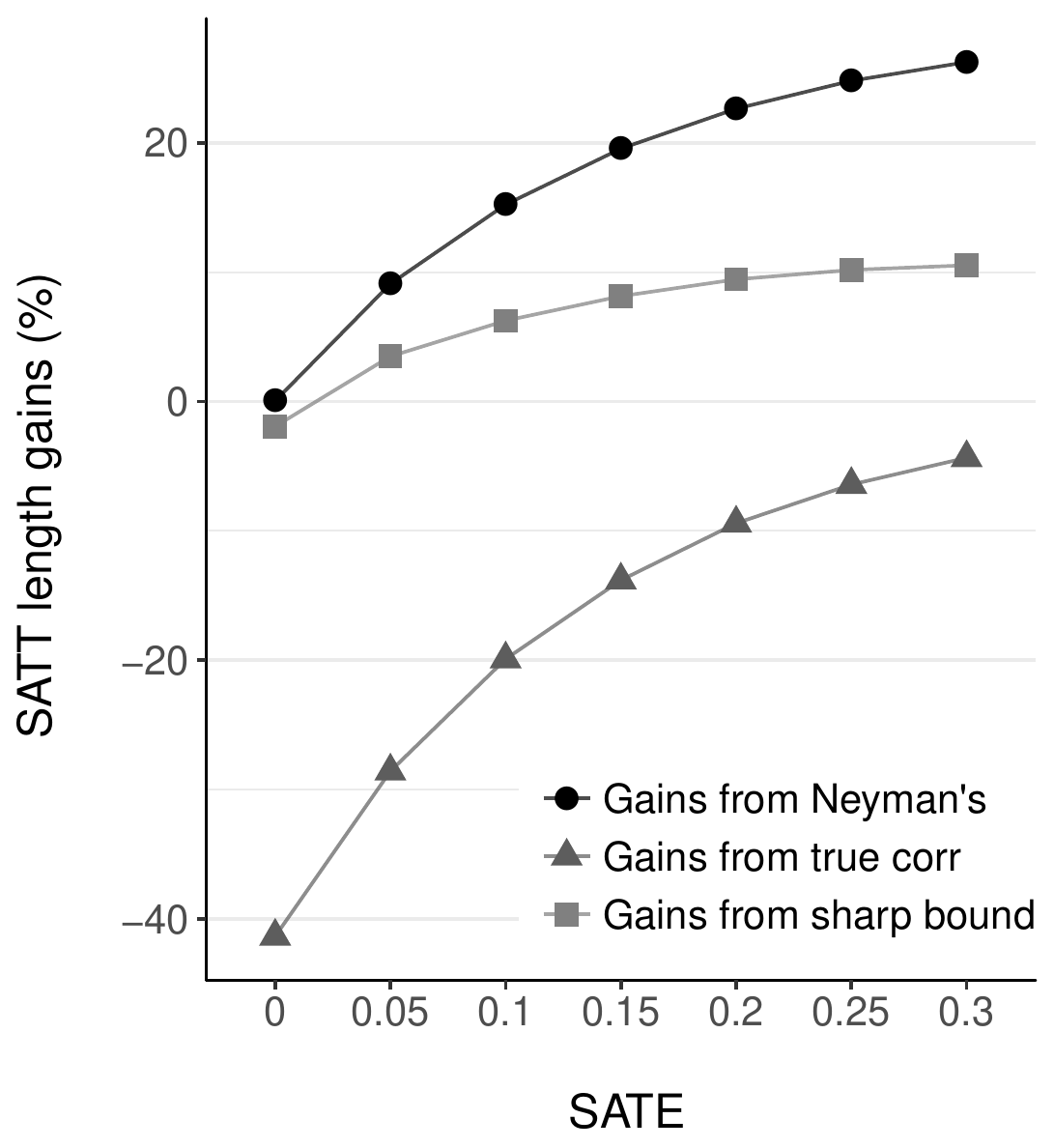}
\begin{minipage}{16cm}
\footnotesize
\emph{Notes}: See the notes in Figure \eqref{fig: random_coefficient simulation}.
\end{minipage}
\end{figure}

Next we decompose how much of the efficiency gain is due to a change of estimand, and how much is due to the conservative variance estimation. The lower line (i.e., triangles) in Figure \eqref{fig: binary outcome gain analytical SATT vs SATE} shows the percentage difference in length relative to a CI for SATE when $\rho$ is known. This CI is substantially shorter than the ones for SATT, which demonstrates that in this simulation the efficiency differences are due to the conservative variance estimation and not because of the change in estimand. Unlike the random coefficient model above, there are small differences in coverage under this data-generating model, although, there are differences in accuracy. The binary outcomes model has less heterogeneity in the treatment effect distribution (i.e., $\sigma_{\tau}^2$ is smaller) which mitigates the differences between the different estimands.       

\subsection{Censored outcomes (Tobit model)}
\label{sec: Analytical example - tobit}

\begin{figure}[ht]
\centering
\caption{ \textbf{Censored outcome (Tobit) simulation results: Length and coverage differences} }
\label{fig: tobit simulation}
\includegraphics[scale=0.55]{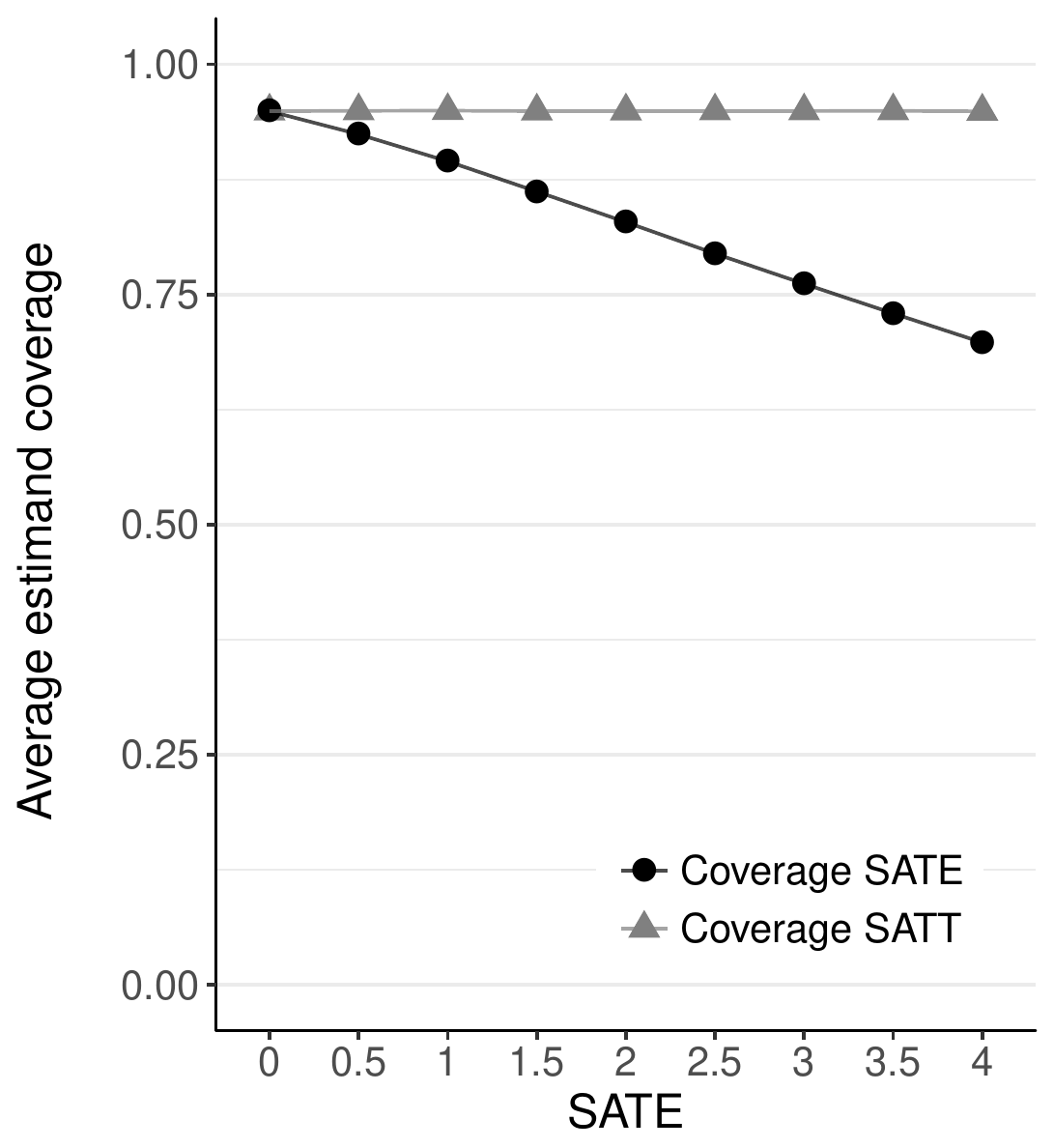}
\includegraphics[scale=0.55]{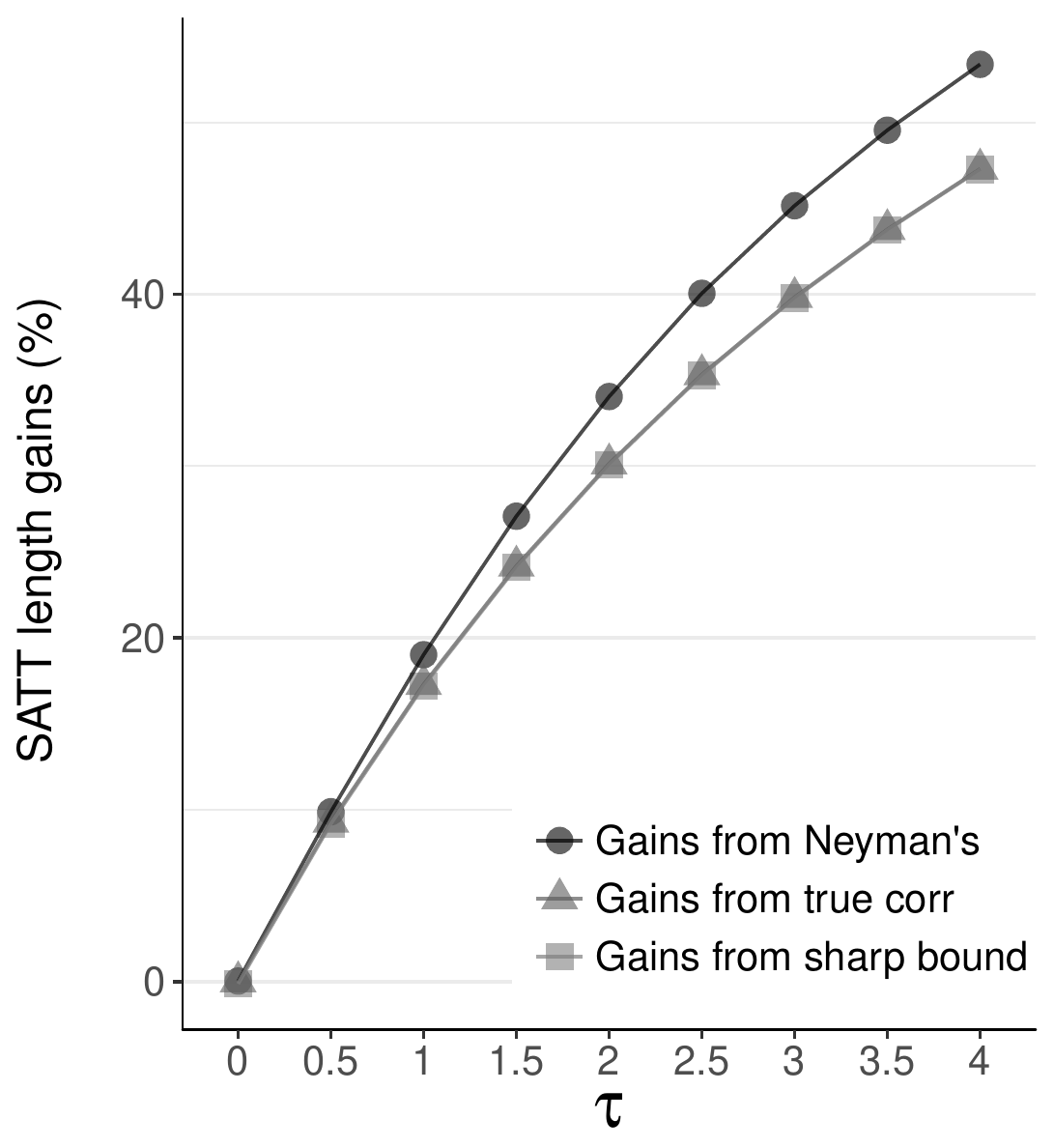}
\begin{minipage}{16cm}
\footnotesize
\emph{Notes}: See the notes in Figure \eqref{fig: random_coefficient simulation}.      
\end{minipage}
\end{figure}

Consider a finite population of size $N$,
\begin{align}
\Pi_N &=  \left\{ ( Y(0)_{1N}, Y(1)_{1N} ), \; 
( Y(0)_{2N}, Y(1)_{2N} ), \dots, \;
( Y(0)_{NN}, Y(1)_{NN} )  \right\}
\end{align}   
where $Y(0)$ is a continuous outcome and $Y(1)$ is,
\begin{align*}
Y(1) = \left\{ 
\begin{array}{ll}
Y(0) + \tau, & Y(0) \geq 0 \\
Y(0), & Y(0) < 0
\end{array} \right. \quad \quad \text{and} \quad \tau > 0
\end{align*}   
The above data generating-process implies that the variance of the potential outcomes under treatment is higher than the variance of units under the control regime:
\begin{align*}
\var{Y(1)} = \sigma_0^2 + \Pr( Y(0) > 0 ) \cdot \tau 
\cdot \left[ \tau \cdot \left( 1 - \Pr( Y(0) > 0 \right) + \e{Y(0)|Y(0) > 0} - \e{Y(0)} \right]
> \sigma_0^2
\end{align*}  
and the variance ratio is:
\begin{align*}
\frac{\sigma_1^2}{\sigma_0^2} = 1 + \frac{ \Pr( Y(0) > 0 ) \cdot \tau 
\cdot \left[ \tau \cdot \left( 1 - \Pr( Y(0) > 0 \right) + \e{Y(0)|Y(0) > 0} - \e{Y(0)} \right] }{ \sigma_0^2 }
\end{align*}
The variance ratio is increasing with respect to $\tau$ and there is a potential for accuracy gains by changing the estimand from SATE to SATT. The simulation results in Figure \eqref{fig: tobit simulation} confirm the above derivation. Similarly to the random coefficient model the change of estimand is the main cause of the differences in length. Similarly to the previous simulations, as $\sigma_{\tau}^2$ increases the differences between estimands become more stark and inference for SATT (SATE) has bad coverage w.r.t. SATE (SATT):     
\begin{align}
\label{eq: tobit dgp}
Y_i(0) & \sim N(\mu = 0, \; \sigma_0^2 = 1) 
\end{align}

\section{Super-population model} \label{sec: Super-population model}

In the super-population model the sampling procedure has two steps. First, a sample of $N$ units is drawn from a super-population, $F_{ Y(1), Y(0) } (\cdot)$, which can be either finite or infinite. Second, $m$ units are randomly allocated to the treatment regime and the remaining $N-m$ units are assigned to the control regime.\footnote{The number of treated and control units is fixed: $m$ is not a random variable.}
Both the treatment indicator and the potential outcomes are random variables, and they are independent, due to the random assignment of units to treatment regimes. SATE, SATT, and SATC are all random variables that are unbiased, and consistent, estimators of PATE. Under the super-population model the difference-in-means test statistic can be used to construct a CI for PATE. Neyman's variance estimator is a consistent estimator for the variance of the difference-in-means test statistic \citep{imbens2015}. PATE, PATT, and PATC are all equal and recentering the difference-in-means w.r.t. PATE is equivalent to recentering it w.r.t. PATT (or PATC). There are no efficiency gains or inaccuracies in making inference on PATT instead of PATE.

Next we describe the behavior of the difference-in-means when it is recentered w.r.t. different sample average treatment effect estimands (e.g., SATE, SATT). The variance of $\td-\text{SATE}$ is lower than that of $\td-\text{PATE}$:
\begin{align}
\label{eq: SATT smaller than PATE variance}
\var{ \td - \text{SATE} } &=
\var{ \td } - \underset{\var{ \text{SATE} }}{\underbrace{\frac{1}{N} \cdot \sigma_{\tau}^2}}  \leq \var{ \td } = \var{ \td - \text{PATE} }
\end{align}
and
\begin{align}
\var{ \td } &=
\frac{\sigma_1^2(\text{SP})}{m} + \frac{\sigma_0^2(\text{SP})}{N-m}  
\end{align}
where $\sigma_j^2(\text{SP})$ ($j = 0,1$) denotes the variance of units under treatment regime $j$ in the super-population. See Appendix \eqref{proof: super-population variances of recentered diff-in-means } for a proof of equation \eqref{eq: SATT smaller than PATE variance}. In the fixed population case, SATE is a parameter and not a random variable, which changes when sampling uncertainty is introduced and the potential outcomes become random variables as well. 

Lemma \eqref{lemma: super-population variances of recentered diff-in-means} establishes that the variances calculations we derived in Section \eqref{sec: Theory fixed/finite population} hold also under the super-population framework. Although the variance calculations are over repeated samples, the correlation term between potential outcomes, $\rho$, does not go away and to conduct inference on the SATE all the issues we discussed in Section \eqref{sec: Theory fixed/finite population} arise here as well. More specifically, the uncertainty calculation over repeated draws of data cancels the $\rho$ term in the variance of the difference-in-means test statistic. However, it also introduces a correlation term $\rho$ in the variance of SATE, which is now a random variable. The two elements cancel each other out and yield a variance formula that is exactly the same as the one that was derived in Lemma \eqref{lemma: variances of recentered diff-in-means}. 
 
\begin{lemma}
\label{lemma: super-population variances of recentered diff-in-means}
Under the super-population model, the variance of the difference-in-means when recentered w.r.t. the SATE, SATT, or SATC is 
\begin{small}
\begin{align*}
\var{ \td - \text{SATE} } &= 
\frac{1}{N \cdot (1-p) \cdot p} \cdot \left[ p^2 \cdot \sigma_0^2(\text{SP}) + (1-p)^2 \cdot \sigma^2_1(\text{SP})
+ 2 p (1-p) \cdot \rho(\text{SP}) \cdot \sigma_0(\text{SP}) \cdot  \sigma_1(\text{SP}) \right]
\\
\var{ \td - \text{SATT} } &= \frac{1}{p \cdot (1-p) \cdot N} \cdot \sigma_0^2(\text{SP}) 
\\
\var{ \td - \text{SATC} } &= \frac{1}{p \cdot (1-p) \cdot N} \cdot \sigma_1^2(\text{SP})
\end{align*}
\end{small}
See Appendix \eqref{proof: super-population variances of recentered diff-in-means } for the proof. 
\end{lemma}
  

\section{Comments} \label{sec: Comments and remark}

Several notes on the previous results and possible extensions are in order.  

\emph{Remark 1}. The previous results can be extended to include covariate adjustment of pre-treatment characteristics according to the procedure that was proposed by \cite{rosenbaum2002}. Denote by $X_i$ a $1 \times p$ dimensional vector of the pre-treatment characteristics of unit $i$. The matrix $\bm{X}$ has dimensions $n \times p$ and each row $i$ contains the pre-treatment characteristics of unit $i$. It is common to adjust $Y_i$ for $X_i$ for efficiency purposes. Define $Y_i^{adjusted} = X_i(X'X)^{-1}X'Y_i$ as the adjusted/residualized responses. All the results for inference on $Y_i$ also apply for inference on $Y_i^{adjusted}$.

\emph{Remark 2}. The inference results for SATT (and SATC) can also be extended to different randomization models. For example, Theorem \ref{theorem: CLT - Bernoulli trials randomization} provides variance and limiting distribution results for inference on SATT when the treatment assignment is done by random independent Bernoulli trials. This illustrates how inference on SATT can be derived under random treatment assignment models that differ from the classic complete randomization mechanism.

\begin{theorem}{(Limiting distribution of \td - SATT)} 
\label{theorem: CLT - Bernoulli trials randomization}
When treatment is assigned according to random independent Bernoulli trials, the standardized and recentered, w.r.t. SATT, difference-in-means follows a standard Normal distribution under two regularity conditions. When the following is satisfied:
\begin{align}
N-m \rightarrow \infty, \quad  m \rightarrow \infty, \quad \text{and} \quad \sigma_1^2, \; \sigma_0^2 < \infty
\end{align}
then:
\begin{align}
\label{eq: CLT - Bernoulli trials randomization}
\frac{ \frac{N-m}{N} \cdot \left( \td - \text{SATT} \right) }{ \sqrt{ \var{ \frac{1}{m} \cdot \sum_{i=1}^N Y_i(0) \cdot T_i } } } \overset{d}{\rightarrow} N(0,1)
\end{align}
where
\begin{align}
\label{eq: variance Bernoulli random assignment}
\var{ \frac{1}{m} \cdot \sum_{i=1}^N Y_i(0) \cdot T_i } 
&=
\left( \begin{array}{c}
\frac{1}{m} \\
 - \frac{\sum_{i=1}^N Y_i(0) \cdot T_i}{m^2}
\end{array} \right)^T
\cdot
\Sigma
\cdot
\left( \begin{array}{c}
\frac{1}{m} \\
 - \frac{\sum_{i=1}^N Y_i(0) \cdot T_i}{m^2}
\end{array} \right)
\end{align}
and
\begin{align}
\Sigma &= \left( 
\begin{array}{cc}
N \cdot p(1-p) \left[ \sigma_0^2 + \e{Y(0)}^2 \right] & p(1-p) \cdot N \cdot  \e{Y(0)} \\
p(1-p) \cdot N \cdot  \e{Y(0)} & N \cdot p(1-p) \\
\end{array}
\right)
\end{align}
See Appendix \eqref{proof: CLT - Bernoulli trials randomization} for the proof. 
\end{theorem}


\emph{Remark 3}. The regularity conditions of the above theoretical results will usually be satisfied in applied research, yet it is important to understand when they will not. For example, imagine a control regime in which all the individuals die ($Y_i(0)=0  \;\; \forall_i$), while under the treatment regime units have a strictly positive survival probability, $\e{Y(1)}=p_1$ and $\var{Y(1)} = p_1 (1-p_1)$. As the variance of the control units is strictly lower (zero) than that of the treated units, there are potential accuracy gains from estimating SATT instead of SATE. In this scenario as $\sigma^2_0 = 0$ the PI for SATT contains only one point, the difference-in-means estimate, which is clearly wrong. The example above does not stand in contradiction to Theorem \eqref{theorem: CLT - fixed sample}, as in the above case the regularity condition of the theorem is not satisfied:
\begin{align*}
\frac{ \underset{ 1 \leq i \leq N }{\max}  \left( Y(0)_{Ni} - \bar{Y}(0)_{N}  \right)^2 }{ \sum_{i=1}^N \left( Y(0)_{Ni} - \bar{Y}(0)_{N}  \right)^2  } = \frac{0}{0}
\end{align*}
and $\frac{0}{0}$ is not a well-defined expression.

\section{Real data application: Online experiments} \label{sec: Main real the application}

To better understand the trade-offs in conducting inference on different average treatment effect estimands, we analyze a sample of online field experiments that have been conducted by a large internet firm as product improvement tests. Our sample consists of 278 experiments with an average sample size of approximately one hundred million units per experiment. Many different outcome metrics are analyzed for various subgroups. The average subgroup consists of 1.1 million observations, and there are 826 unique outcome metrics across all of the experiments. In total, twenty-five thousand different treatment effects are estimated. The data analyzed were aggregated and de-identified.

\begin{figure}[ht]
\centering
\caption{ \textbf{ Distribution of the difference-in-means across experiments } }
\label{fig: raw difference in means distribution}
\includegraphics[scale=0.75]{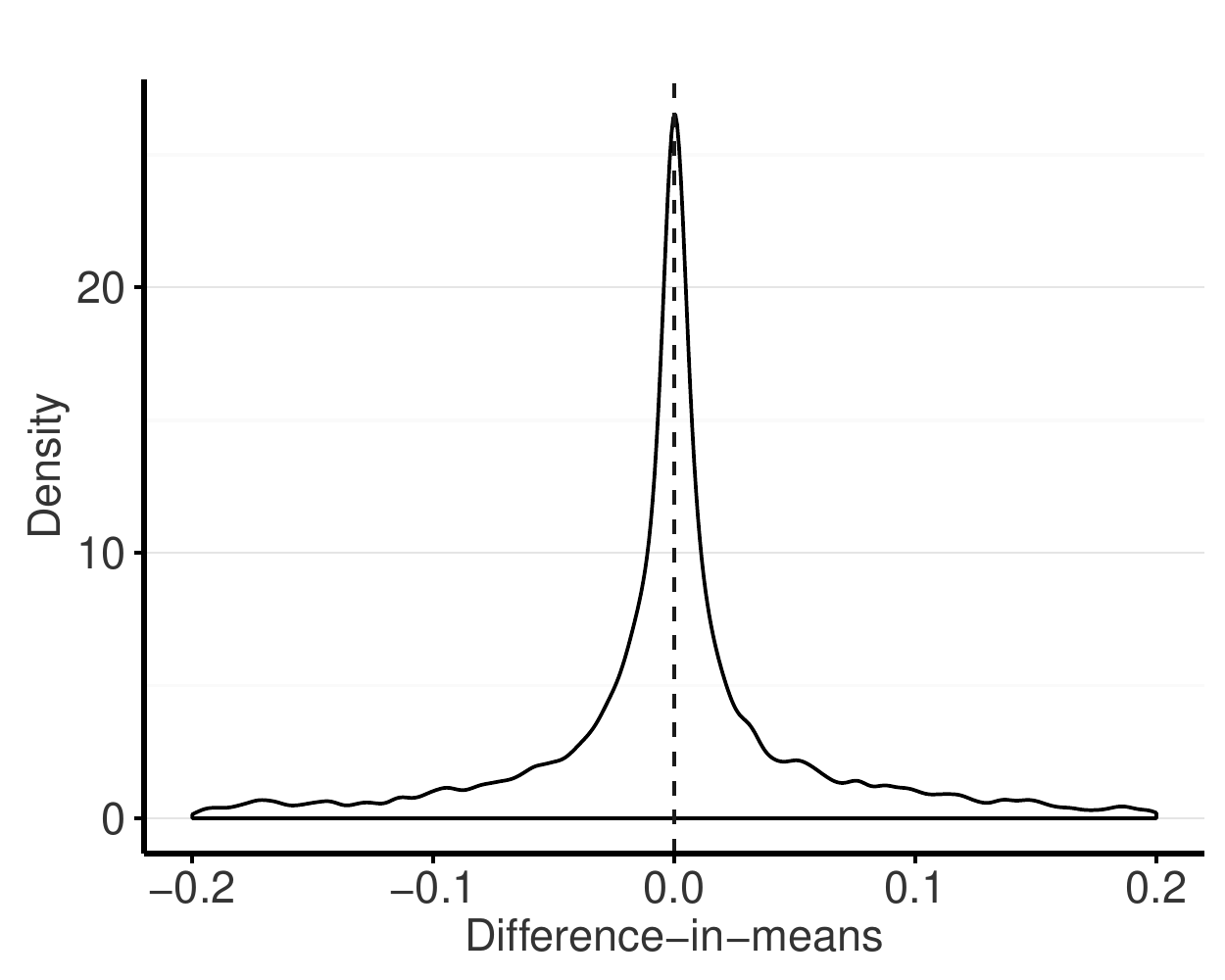}
\begin{minipage}{10cm}
\footnotesize
\emph{Notes}: The figure shows the distribution of the raw difference-in-means over all online experiments in our sample.  
\end{minipage}
\end{figure}

Figure \eqref{fig: raw difference in means distribution} shows the distribution of the difference-in-means across all of the online experiments. It is clear that on average, across experiments, the treatment had a zero effect. The distribution of difference-in-means is tightly centered around zero. However, this does not imply that the treatment had no effect. Treatment effect heterogeneity can generate positive effects for some units and negative effects for others that cancel each other on average. Elaborate tests and computations are difficult to carried out with millions of observations. Inference on SATO provides a simple solution as it is sensitive to heterogeneous treatment effects and is guaranteed to have correct coverage and valid inference on SATO. Another advantage of SATO is that for the simple cases of $\rho \in \lbrace -1,0,1 \rbrace$, the change of estimand from SATE to SATO does not require any additional information for inference. It can be carry out using only aggregate summary statistics that are already collected for inference on SATE.

\begin{figure}[ht]
\centering
\caption{ \textbf{Inference accuracy: A comparison of SATO and SATT to SATE across online experiments}  }
\label{fig: internet experiments - accuracy}
    \begin{subfigure}[b]{0.45\textwidth}
        {\includegraphics[width=\textwidth]{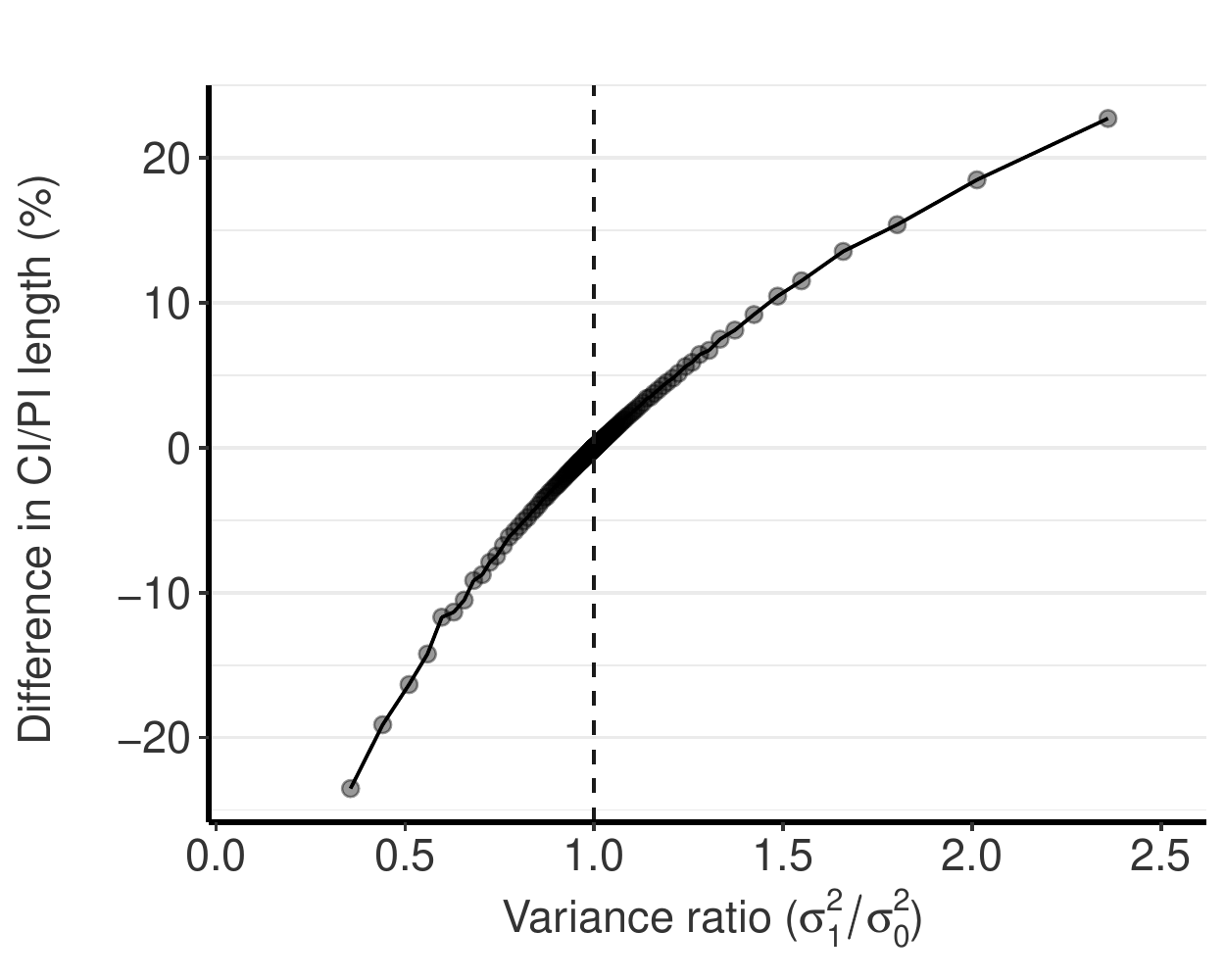}}
        \caption{  Differences in variance: \\ SATT vs. SATE }
    \end{subfigure}
    \begin{subfigure}[b]{0.45\textwidth}
        {\includegraphics[width=\textwidth]{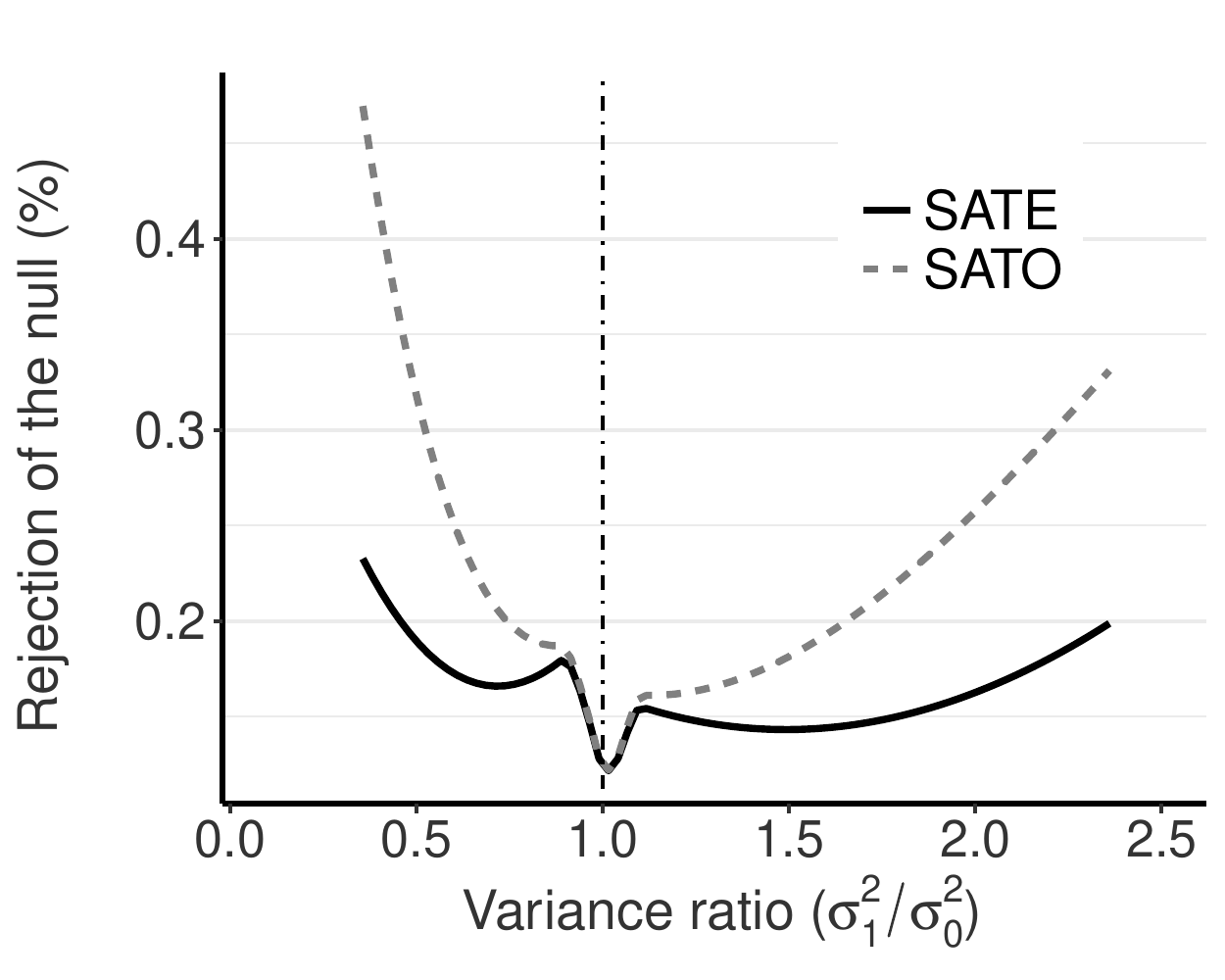}}
        \caption{  Differences in rejection rate: \\  SATO vs. SATE }
    \end{subfigure}
\begin{minipage}{16cm}
\footnotesize
\emph{Notes}: The left plot describes the empirical analogue of Figure \eqref{fig: gain analytical SATT vs SATE}. It plots the relationship between the variance ratio in the experiment and the gain in accuracy from conducting inference on SATT relative to SATE. The right plot compares the rejection rate of the null hypothesis that the estimand of interest is equal to zero for SATE and SATO when $\rho = 1$.         
\end{minipage}
\end{figure}

The left plot of Figure \eqref{fig: internet experiments - accuracy}  replicates Figure \eqref{fig: gain analytical SATT vs SATE} in a real empirical application and demonstrates that variance differences across experiments exist. As the variance ratio $\frac{\sigma_1^2}{\sigma_0^2}$ increases, the PI for SATT becomes shorter, and in some cases there can be large variance gains from changing the estimand to SATO. The left plot of Figure \eqref{fig: internet experiments - accuracy} shows the share of experiments in which the null hypothesis that the estimand is equal to zero was rejected, as a function of the variance ratio. As the variance ratio departs from one there are substantial efficiency gains to be had from conducting inference on SATO (with $\rho=1$) relative to SATE. The differences in rejection rates between inference on SATE and SATO are evidence of treatment effect heterogeneity.

\FloatBarrier
\section{Discussion} \label{sec: Discussion} 

In some applications, the estimand of interest is SATT (or SATC). Making inferences about SATT (or SATC) using a CI for SATE relies on variance estimators that are not consistent and are not guaranteed to have correct coverage or be efficient. We derive efficient variance formulas for inference on a new and general class of estimands derived from any mixing between SATT and SATC. The variance formulas are used to construct PIs that are non-parametrically guaranteed to have correct coverage and to be non-conservative -- unlike inference on SATE. All inference procedures discussed in the paper use the difference-in-means as the test statistic, and therefore have the same point estimates as existing methods. Note that all three estimands, SATE, SATT, and SATC, are equal in expectation. The key difference is in the variance calculations. In addition, we present a new diagnostic tool, a decomposition of the uncertainty about the value of SATE into two components: uncertainty about SATT and about SATC, which both have an intuitive interpretation.

Taken together, the Monte Carlo simulations demonstrate that: (i) the choice of estimand has a direct implication on the accuracy of the inference that can be conducted; (ii) using variance formulas that have correct size for SATE will \emph{not} have correct coverage of other sample average treatment effect estimands such as SATT or SATO. When a researcher is interested in conducting inference on SATT, she is required to use variance formulas that have correct coverage for this estimand. Such formulas are derived in the paper.    

The application to online experiments is a clear example of a scenario in which the methods developed in the paper can provide accuracy gains over standard practice. The gains mainly arise because inference on SATO is sensitive to treatment effect heterogeneity through the variance calculations. The large potential for accuracy gains emphasizes that researchers and other applied users should think carefully about which causal estimand they want to conduct inference on. If the answer is SATE, then they should not use any of the results in this paper; however, if SATO is also of interest, then large accuracy gains can be obtained by a change of estimand. In addition, if SATT is the primary object of interest, then variance formulas for inference on SATE will have incorrect coverage of SATT.     

Our results can be used in a variety of empirical applications from medicine and social sciences to online ad experiments.  In many cases one would like to have more accuracy when testing the null hypothesis of no average treatment effect. This occurs even when, as in the case of digital experiments, the sample size is large, but the power is limited \citep{lewis2014unfavorable}. The accuracy gains we document hold even asymptotically, and could be of special value in massive experiments because our approach  requires computing only summary statistics that are simple to obtain and easy to compute, and lead to possibly large gains in accuracy.

\clearpage
\singlespacing
\bibliography{Bibstat}
\bibliographystyle{aer}
\clearpage

\section*{\LARGE Online appendix: Proofs}
\appendix
\doublespacing

\section{Derivation of MSE(SATE,SATT)}
\label{proof: derivation of MSE(SATE,SATT)}

\begin{align*}
\text{MSE}(\text{SATE}, \; \text{SATT})  &= 
\text{MSE}( \td -  \text{SATE} , \; \td -  \text{SATT}) 
\\ &=
\left[ \e{ \text{SATE} - \text{SATT} } \right]^2
+
\var{ \text{SATT} } 
\nonumber \\ &= 
\var{ \text{SATT} } 
\nonumber \\ &= 
\left( \frac{1}{m} - \frac{1}{N} \right) \cdot \sigma^2_{\tau}
\nonumber \\ &= 
\frac{N-m}{Nm} \cdot \var{\tau_i} = \frac{1-p}{m} \cdot \sigma^2_{\tau}
\end{align*} 

\section{Proof of Lemma \ref{lemma: decomposition of the difference in means}}
\label{proof: diff-in-means decomposition}

The sum of the responses of the treated units can be decomposed into two components, the sum of the potential outcomes under control and the sum of the treatment effects. 
\begin{align}
\label{eq: decomposition of treated sum}
\sum_{i=1}^N Y_i \cdot T_i 
&= 
\sum_{i=1}^N T_i \cdot Y_i(1) 
= \sum_{i=1}^N T_i \cdot Y_i(0) + \sum_{i=1}^N T_i \cdot (Y_i(1) - Y_i(0))
\\ &\Rightarrow 
\frac{1}{m} \cdot \sum_{i=1}^N Y_i \cdot T_i 
= 
\frac{1}{m}  \cdot \sum_{i=1}^N T_i \cdot Y_i(0) + \text{SATT} 
\end{align}
We can re-write the difference in means as,
\begin{align}
\td &= \frac{1}{m} \cdot \sum_{i=1}^N Y_i \cdot T_i - \frac{1}{N-m} \cdot \sum_{i=1}^N Y_i \cdot (1-T_i)  \nonumber \\
&= \frac{1}{m}  \cdot \sum_{i=1}^N T_i \cdot Y_i(0) + \text{SATT} - \frac{1}{N-m} \cdot \sum_{i=1}^N Y_i(0) \cdot (1-T_i)
\nonumber \\ &= 
\frac{1}{m}  \cdot \sum_{i=1}^N T_i \cdot Y_i(0) - \frac{1}{N-m} \cdot \sum_{i=1}^N Y_i(0) \cdot (1-T_i) + \text{SATT}
\nonumber \\
&= \frac{N}{m \cdot (N-m)} \cdot \sum_{i=1}^N Y_i(0) T_i -\frac{1}{N-m} \cdot \underset{\text{constant}}{\underbrace{\sum_{i=1}^N Y_i(0)}} + \text{SATT}
\end{align}

The decomposition of $\td$ as a function of $Y(1)$ and the SATC can be proved using a similar proof to the one above. 

\section{ Proof of Lemma \ref{lemma: variances of recentered diff-in-means} }
\label{proof: variance of different recentered diff-in-means}

The variance of $\td - \text{SATE}$ can be written as,

\begin{align*}
\var{ \td - \text{SATE} } &=  \frac{m}{N \cdot (N-m)} \sigma_0^2 + \frac{N-m}{N \cdot m} \sigma_1^2 
+ \frac{2}{N} \cdot \rho \cdot \sigma_0 \cdot  \sigma_1
\\ &=
\frac{p}{N \cdot (1-p)} \cdot \sigma_0^2 + \frac{N \cdot (1-p)}{N^2 \cdot p} \cdot \sigma^2_1
+ \frac{2}{N} \cdot \rho \cdot \sigma_0 \cdot  \sigma_1
\\ &=
\frac{p}{N \cdot (1-p)} \cdot \sigma_0^2 + \frac{ (1-p)}{N \cdot p} \cdot \sigma^2_1
+ \frac{2}{N} \cdot \rho \cdot \sigma_0 \cdot  \sigma_1
\\ &=
\frac{1}{N \cdot (1-p) \cdot p} \cdot \left[ p^2 \cdot \sigma_0^2 + (1-p)^2 \cdot \sigma^2_1
+ 2 p (1-p) \cdot \rho \cdot \sigma_0 \cdot  \sigma_1 \right]
\end{align*}  
The variance derivations of $\td - \text{SATT}$ are described in \ref{proof: CLT - fixed sample}, and the variance calculation of $\td - \text{SATC}$ follows exactly the same steps as the variance calculation of $\td - \text{SATT}$.

\section{ Proof of Theorem \ref{theorem: gains from estimand} } 
\label{proof: gains from estimand change}

\subsubsection*{Proof of part (1) of Theorem \ref{theorem: gains from estimand} }

The proof is a simple implementation of the intermediate value theorem. We show that the difference in variance is negative when $\rho = -1$, and is positive when $\rho = 1$, and that the difference in variance is a continuous and increasing function with respect to $\rho$. Then the desired result follows immediately from the intermediate value theorem. 

When $\rho = -1$,  the variance of $\td - \text{SATE}$ can be written as:
\begin{align*}
\var{ \td - \text{SATE} } &=  \frac{m}{N \cdot (N-m)} \sigma_0^2 + \frac{N-m}{N \cdot m} \sigma_1^2 
+ \frac{2}{N} \cdot \rho \cdot \sigma_0 \cdot  \sigma_1
\\ &=
\frac{p}{N \cdot (1-p)} \cdot \sigma_0^2 + \frac{N \cdot (1-p)}{N^2 \cdot p} \cdot \sigma^2_1
- \frac{2}{N} \cdot \sigma_0 \cdot  \sigma_1
\\ &=
\frac{p}{N \cdot (1-p)} \cdot \sigma_0^2 + \frac{ (1-p)}{N \cdot p} \cdot \sigma^2_1
- \frac{2}{N} \cdot \sigma_0 \cdot  \sigma_1
\\ &=
\frac{1}{N \cdot (1-p) \cdot p} \cdot \left[ p^2 \cdot \sigma_0^2 + (1-p)^2 \cdot \sigma^2_1
- 2 p (1-p) \cdot \sigma_0 \cdot  \sigma_1 \right]
\\ &= 
\frac{1}{N \cdot (1-p) \cdot p} \cdot \left( p \cdot \sigma_0 - (1-p) \cdot \sigma_1 \right)^2
\end{align*}  
and as, 
\begin{align*}
p \cdot \sigma_0 < \sigma_0 - (1-p) \cdot \sigma_1 < \sigma_0
\end{align*}
it follows that when $\rho = -1$, then $\var{ \td - \text{SATE} } < \var{ \td - \text{SATT} }$.  
When $\rho = 1$,  the variance of $\td - \text{SATE}$ can be written as:
\begin{align*}
\var{ \td - \text{SATE} } &=  \frac{p}{N \cdot (1-p)} \cdot \sigma_0^2 + \frac{N \cdot (1-p)}{N^2 \cdot p} \cdot \sigma^2_1
+ \frac{2}{N} \cdot \sigma_0 \cdot  \sigma_1
\\ &=
\frac{p}{N \cdot (1-p)} \cdot \sigma_0^2 + \frac{ (1-p)}{N \cdot p} \cdot \sigma^2_1
+ \frac{2}{N} \cdot \sigma_0 \cdot  \sigma_1
\\ &=
\frac{1}{N \cdot (1-p) \cdot p} \cdot \left[ p^2 \cdot \sigma_0^2 + (1-p)^2 \cdot \sigma^2_1
+ 2 p (1-p) \cdot \sigma_0 \cdot  \sigma_1 \right]
\\ &= 
\frac{1}{N \cdot (1-p) \cdot p} \cdot \left( p \cdot \sigma_0 + (1-p) \cdot \sigma_1 \right)^2
\end{align*}
and as $\sigma_1 > \sigma_0$, it follows immediately that when $\rho = 1$, $\var{ \td - \text{SATE} } > \var{ \td - \text{SATT} }$.  

As $\frac{\partial \var{\td - \text{SATE}} - \var{\td - \text{SATT}}}{ \partial \rho } > 0$, it follows directly from the intermediate value theorem that there exists a value of $\rho$, $\bar{\rho}$, such that:
\begin{align*}
\rho \leq  \bar{\rho} \Rightarrow \var{ \td - \text{SATE} } \leq \var{ \td - \text{SATT} }
\\
\rho >  \bar{\rho} \Rightarrow \var{ \td - \text{SATE} } > \var{ \td - \text{SATT} }
\end{align*}

\subsubsection*{Proof of part (2) of Theorem \ref{theorem: gains from estimand} }

Consider the case in which $\rho = 0$: 
\begin{align*}
\var{ \td - \text{SATE} } &= 
\frac{1}{N \cdot (1-p) \cdot p} \cdot \left[ p^2 \cdot \sigma_0^2 + (1-p)^2 \cdot \sigma^2_1 \right]
\end{align*}  
$\var{ \td - \text{SATT} }$ will be lower than $\var{ \td - \text{SATE} }$, when:
\begin{align*}
p^2 \cdot \sigma_0^2 + (1-p)^2 \cdot \sigma^2_1 \geq \sigma_0^2
\Rightarrow
\frac{(1-p)^2}{1-p^2} \cdot \sigma^2_1 \geq \sigma_0^2
\end{align*}
Note, $\frac{(1-p)^2}{1-p^2} = \frac{1+p^2 -2p}{1-p^2} \leq 1 
\Rightarrow 
p^2 - 2p \leq -p^2
\Rightarrow 
2p^2 \leq 2p
\Rightarrow 
p^2 \leq p
$, which is always satisfied. Therefore, $ 0 \leq \frac{(1-p)^2}{1-p^2} \leq 1 $. To conclude, when $\frac{\sigma_1}{\sigma_0} > \sqrt{\frac{1-p^2}{(1-p)^2}}$  the variance difference will be positiv $\var{ \td - \text{SATE} } > \var{ \td - \text{SATT} }$. As the variance difference is negative when $\rho = -1$ it follows directly from the intermediate value theorem that the desired $\bar{\rho}$ exists and is strictly larger than $-1$ and lower than $0$ --- $\bar{\rho}$ is negative. This concludes the proof.

\section{Proof of Theorem \ref{theorem: prediction vs. confidence intervals - fixed sample} }
\label{proof: CLT - fixed sample}
According to Lemma \ref{lemma: decomposition of the difference in means}, the adjusted difference in means, $\frac{N-m}{N} \cdot \left( \td - \att \right)$, has the same distribution of $\frac{1}{m} \cdot \sum_{i=1}^N Y_i(0) \cdot T_i$ with an additive shift of $\frac{1}{N} \cdot \sum_{i=1}^N Y_i(0)$. \cite{hajek1960}, \cite{hajek1961}, \cite{lehmann1975}, and \cite{peng2016} all provide proofs showing that in a finite population model with complete randomization the standardized mean of the treated (or sampled) units will follow a standard Normal distribution under two regularity conditions, 
\begin{align}
N-m \rightarrow \infty, \quad \text{and} \quad m \rightarrow \infty
\end{align}
and 
\begin{align}
\frac{ \underset{ 1 \leq i \leq N }{\max}  \left( Y(0)_{Ni} - \bar{Y}(0)_{N}  \right)^2 }{ \sum_{i=1}^N \left( Y(0)_{Ni} - \bar{Y}(0)_{N}  \right)^2  }
\cdot
\max \left( \frac{N-m}{m}, \frac{m}{N-m}  \right)
\rightarrow \infty
\end{align}

The expectation and variance of $\frac{N-m}{N} \cdot \left( \td - \att \right)$ are derived next and complete the proof of (\ref{theorem: CLT - fixed sample}). The expectation is,
\begin{align*}
\e{\frac{N-m}{N} \cdot \left( \td - \att \right)} 
&=
\frac{N-m}{N} \cdot \left[ \frac{N}{m(N-m)} \cdot \sum_{i=1}^N \e{Y(0)} \cdot \e{T} - \frac{1}{N} \cdot \sum_{i=1}^N \e{Y(0)} \right]
\\ &=
\frac{N}{m} \cdot \e{Y(0)} \cdot \frac{m}{N} - - \frac{1}{N} \cdot N \cdot \e{Y(0)}
\\ &= 
0
\end{align*} 
where the first equality follows from Lemma \ref{lemma: decomposition of the difference in means}. The variance of the adjusted difference in means is,
\begin{align*}
\var{ \frac{N-m}{N} \cdot \left( \td - \att \right) } 
&= 
\var{ \frac{1}{m} \cdot \sum_{i=1}^N Y_i(0) \cdot T_i - \frac{1}{N} \cdot \sum_{i=1}^N Y_i(0) } 
\\ &=
\frac{1}{m^2} \cdot \var{ \cdot \sum_{i=1}^N Y_i(0) \cdot T_i}
\end{align*} 
where the first equality follows from Lemma \ref{lemma: decomposition of the difference in means}, and the second equality follows as $\frac{1}{N} \cdot \sum_{i=1}^N Y_i(0)$ is a constant and not a random variable. The variance of the sum $\sum_{i=1}^N Y_i(0) \cdot T_i$, is the variance of $m$ units who are randomly sampled from the population $Y_{N1}(0), \dots, Y_{NN}(0)$.\footnote{Our calculation of the variance of $\frac{1}{N} \cdot \sum_{i=1}^N Y_i(0)$ follows the variance calculation of the Wilcoxon Rank Sum Test statistic in \cite{lehmann1975}.} Denote the $m$ chosen variables by $V_1, \dots, V_m$, and the variance of the sum $V_1+ \dots + V_m$ is,
\begin{align}
\label{eq: variance of sum Y(0) * T}
\var{ V_1 + \dots + V_m } 
&= 
\sum_{i=1}^m \var{V_i} + \sum_{i=1}^m \sum_{j \neq i} \cov{V_i,V_j}
= m \cdot \sigma_0^2 + m(m-1) \cdot \rho
\end{align} 
where the second equality follows as the covariance between each two units is the same and is denoted by $\rho$. When $m=N$, the variance of the sum is zero and therefore,
\begin{align}
\label{eq: covariance between Y(0) of controls}
\rho &= - \frac{\sigma_0^2}{N-1}
\end{align} 
\noindent Substituting \eqref{eq: covariance between Y(0) of controls} in \eqref{eq: variance of sum Y(0) * T} and re-arranging yields that:
\begin{align*}
\var{ \frac{N-m}{N} \cdot \left( \td - \att \right) } 
&= 
\var{ \frac{1}{m} \cdot \sum_{i=1}^N Y_i(0) \cdot T_i - \frac{1}{N} \cdot \sum_{i=1}^N Y_i(0) } 
\\ &=
\frac{1}{m^2} \cdot \var{ \cdot \sum_{i=1}^N Y_i(0) \cdot T_i} 
=
\frac{1}{m^2} \cdot \frac{m(N-m)}{N-1} \cdot \sigma_0^2
\\ &=
 \frac{N-m}{m(N-1)} \cdot \sigma_0^2
\\ & \approx \frac{N-m}{m N} \cdot \sigma_0^2
=
\left( \frac{1}{m} - \frac{1}{N} \right) \cdot \sigma_0^2
\end{align*}
Hence, it follows that:
\begin{align*}
\var{ \left( \td - \att \right) } 
&= \left( \frac{N}{N-m} \right)^2 \cdot \var{ \frac{N-m}{N} \cdot \left( \td - \att \right) } 
\\ &=
\left( \frac{N}{N-m} \right)^2 \cdot \left( \frac{N-m}{mN} \right) \cdot \sigma_0^2
\\ &= 
\frac{N}{(N-m) \cdot m} \cdot \sigma_0^2
\\ &= \frac{1}{N p (1-p)} \cdot \sigma_0^2
\end{align*}

\section{Proof of Theorem \eqref{theorem: Combining PIs for CI for SATE - corr=1} }
\label{proof: Combining PIs for CI for SATE - corr=1}

The first part of the theorem, that $\left[  p \cdot L_{\att}  + (1-p) \cdot L_{\atc}, \quad p \cdot U_{\att}  + (1-p) \cdot U_{\atc}  \right]$, follows immediately from the proof of Theorem (1) in \cite{rigdon2015}. Next we show the derivations of the second part of the theorem.
\begin{align*}
p \cdot L_{\att}  + (1-p) \cdot L_{\atc} 
&=
\td - z_{1-\alpha/2} \cdot \sqrt{k(N,m)} \cdot \left( p \sigma_0 + (1-p) \sigma_1 \right)  
\\ &=
\td - z_{1-\alpha/2} \cdot \sqrt{k(N,m)} \cdot \sqrt{\hat{\mathbb{V}}_{\rho=1}}  
\end{align*}   
and
\begin{align*}
p \cdot U_{\att}  + (1-p) \cdot U_{\atc} 
&=
\td + z_{1-\alpha/2} \cdot \sqrt{k(N,m)} \cdot \left( p \sigma_0 + (1-p) \sigma_1 \right)  
\\ &=
\td + z_{1-\alpha/2} \cdot \sqrt{k(N,m)} \cdot \sqrt{\hat{\mathbb{V}}_{\rho=1}}  
\end{align*}   
The second equality in both of the above equations follows from substituting $\rho = 1$ in the variance formula for $\var{ \td - \ate }$  in Lemma \eqref{lemma: variances of recentered diff-in-means}. 

\section{Proof of Theorem \ref{theorem: prediction vs. confidence intervals - fixed sample} }
\label{proof: prediction vs. confidence intervals - fixed sample}

The desired result follows immediately from a comparison of the different variance expressions in Lemma \ref{lemma: variances of recentered diff-in-means}.

\section{Super-population model}
\label{proof: super-population results}

\subsection{Proof of Lemma \ref{lemma: super-population variances of recentered diff-in-means} }
\label{proof: super-population variances of recentered diff-in-means }

\begin{align*}
\var{ \td - \text{PATE} } &= \var{ \td } = \frac{\sigma_0^2(\text{SP})}{N-m} + \frac{\sigma_1^2(\text{SP})}{m}
\end{align*}
see \cite{imbens2015} for proof. 

Next we derive the variance of $\td - \text{SATE}$:
\begin{align*}
\var{ \td - \text{SATE} } &= \e{ \left( \td - \text{SATE} \right)^2 } - \left( \td - \text{SATE} \right)^2 
\\ &= \
\e{ \left( \td - \text{SATE} \right)^2 }
\\ &= 
\e{ \td^2 } + \e{ \text{SATE}^2 } - 2 \cdot \e{\td \cdot \text{SATE} }
\\ &= 
\e{ \td^2 } + \e{ \text{SATE}^2 } - 2 \cdot \left( \var{\text{SATE}} + \text{PATE}^2  \right)
\\ &= 
\var{ \td } + \var{ \text{SATE} } - 2 \cdot \var{\text{SATE}}
\\ &= 
\var{ \td } - \var{ \text{SATE} }
\end{align*}
Note that, 
\begin{align*}
\e{\td \cdot \text{SATE} } &= \mathbb{E}_{\text{SP}} \left[ \mathbb{E}_{\text{T}} \left[ \td \cdot \text{SATE} | \left( \mathbf{Y(1)}, \mathbf{Y(0)} \right) \right]  \right] 
\\&= 
\mathbb{E}_{\text{SP}} \left[ \text{SATE}^2 \right] 
\\&=  
\var{\text{SATE}} + \left[ \mathbb{E}_{\text{SP}} \left[ \text{SATE} \right]  \right]^2
\\&=
\var{\text{SATE}} + \text{PATE}^2 
\end{align*}

The variance of SATE is:
\begin{align*}
\var{\text{SATE}} = \frac{1}{N} \cdot \var{Y_i(1) - Y_i(0)} 
= 
\frac{1}{N} \cdot \left( \sigma_1^2(\text{SP}) + \sigma_0^2(\text{SP}) - 2 \rho(\text{SP}) \sigma_1(\text{SP}) \sigma_0(\text{SP})  \right)
\end{align*}

Therefore the $\var{ \td - \text{SATE} }$ is:
\begin{small}
\begin{align*}
\var{ \td - \text{SATE} } 
&=
\frac{1}{m} \cdot \sigma_1^2(\text{SP}) + \frac{1}{N-m} \cdot \sigma_0^2(\text{SP}) - \frac{1}{N} \cdot \left( \sigma_1^2(\text{SP}) + \sigma_0^2(\text{SP}) - 2 \rho(\text{SP}) \sigma_1(\text{SP}) \sigma_0(\text{SP})  \right)
\\ &=
\frac{N-m}{Nm} \cdot \sigma_1^2(\text{SP}) + \frac{m}{(N-m)N} \cdot \sigma_0^2(\text{SP}) + \frac{2 \rho(\text{SP}) \sigma_1(\text{SP}) \sigma_0(\text{SP})}{N}
\\ &=
\frac{1-p}{Np} \cdot \sigma_1^2(\text{SP})  + \frac{p}{N(1-p)} \cdot \sigma_0^2(\text{SP}) + \frac{2 \rho(\text{SP}) \sigma_1(\text{SP}) \sigma_0(\text{SP})}{N}
\\ &=
\frac{1}{Np(1-p)} \cdot \left[ (1-p)^2 \cdot \sigma_1^2(\text{SP}) + p^2 \cdot \sigma_0^2(\text{SP}) \right]
+ \frac{2\rho(\text{SP}) \sigma_1(\text{SP}) \sigma_0(\text{SP})}{N}
\\ &=
\frac{1}{N \cdot (1-p) \cdot p} \cdot \left[ p^2 \cdot \sigma_0^2(\text{SP}) + (1-p)^2 \cdot \sigma^2_1(\text{SP})
+ 2 p (1-p) \cdot \rho(\text{SP}) \cdot \sigma_0(\text{SP}) \cdot  \sigma_1(\text{SP}) \right]
\end{align*}
\end{small}

Next we derive the variance of $\td - \text{SATT}$ under the super-population model. 

\begin{align*}
\var{ \frac{N-m}{N} \cdot \left( \td - \text{SATT} \right) } &= 
\ee{ \vv{ \frac{N-m}{N} \cdot \left( \td - \text{SATT} \right)  | Y(1),Y(0) }{T} }{Y(1),Y(0)}
\\ & \quad +
\vv{ \ee{ \frac{N-m}{N} \cdot \left( \td - \text{SATT} \right)  | Y(1),Y(0) }{T} }{Y(1),Y(0)}
\\ &= 
\ee{ \vv{ \frac{N-m}{N} \cdot \left( \td - \text{SATT} \right)  | Y(1),Y(0) }{T} }{Y(1),Y(0)}
\\ &= 
\ee{ \frac{N-m}{m(N-1)} \cdot \sigma_0^2(\text{FS}) }{Y(1),Y(0)}
\\ &= 
\frac{N-m}{m(N-1)} \cdot \ee{ \sigma_0^2(\text{FS}) }{Y(1),Y(0)}
\\ &= 
\frac{N-m}{m(N-1)} \cdot \sigma_0^2(\text{SP})
\end{align*}
where the first equality follows from the law of total variance (i.e., law of conditional variance). There are a few more technical steps to complete the derivation. We neglect them as they already appeared in detail in Appendix \ref{proof: variance of different recentered diff-in-means}. The calculation of the variance of $\td - \text{SATC}$ is similar to the one above.

\section{Proof: Variance of $( \td - \text{SATO} )$}
\label{proof: diff-in-means recentered w.r.t SATO variance}

Consider the following class of average treatment effect estimands:
\begin{align}
\label{eq: SATO definition}
\ato = \omega \cdot \att + (1-\omega) \cdot \atc
\end{align} 
We discuss inferences on this class of estimands using the difference-in-means test statistic. 
\begin{align*}
\var{ \td - \ato } &=  \var{\td} + \var{\ato} - 2 \cdot \cov{ \ato, \; \td }
\end{align*} 
According to the derivations that are detailed later on in this section we can re-write the variance of $\td$ re-centered w.r.t SATO as:
%
%
\begin{align*}
\var{ \td - \ato } &=  \frac{1}{Np(1-p)} \left[  p^2 \sigma_0^2 + (1-p)^2 \sigma_1^2 + 2 p(1-p) \cdot \rho \sigma_0 \sigma_1 \right]
\\ &\quad 
+ \frac{\sigma_{\tau}^2}{N-1} \cdot \left[ \omega^2 \cdot \frac{1-p}{p} + (1-\omega)^2 \cdot \frac{p}{1-p} -2 \omega (1-\omega)  \right]
\\ &\quad  
- 2 \cdot \left[ \frac{\omega}{(N-1)p} \left[ \sigma_1^2-\rho \sigma_1 \sigma_0   \right] - \frac{1}{N-1} \cdot \sigma_{\tau}^2 + \frac{(1-\omega)}{ (N-1)(1-p) } \left[ \sigma_0^2 - \rho \sigma_1 \sigma_0  \right]  \right]
\end{align*}

\subsection*{Calculation of $\var{ \ato }$ }

\begin{align*}
\var{\ato} &= \omega^2 \cdot \var{\att} + (1-\omega)^2 \cdot \var{\atc}
- 2 \cdot \omega (1-\omega) \cdot \cov{\att,\; \atc}
\\ &=
\omega^2 \cdot \frac{N-m}{m(N-1)} \cdot \sigma_{\tau}^2 
+
(1-\omega)^2 \cdot \frac{m}{(N-m)(N-1)} \cdot \sigma_{\tau}^2 
- 
\omega (1-\omega) \cdot \frac{2}{N-1} \cdot \sigma_{\tau}^2
\\ &=
\frac{\sigma_{\tau}^2}{N-1} \cdot \left[ \omega^2 \cdot \frac{1-p}{p} + (1-\omega)^2 \cdot \frac{p}{1-p} -2 \omega (1-\omega) \cdot \right]
\end{align*}

\subsection*{Calculation of $\cov{ \ato, \; \td }$}
The covariance term is:
\begin{align*}
& \cov{ \omega \cdot \att + (1-\omega) \cdot \atc, \; \frac{N}{m \cdot (N-m)} \cdot \sum_{i=1}^N Y_i(0) \cdot T_i -\frac{1}{N-m} \cdot \sum_{i=1}^N Y_i(0) + \text{SATT} }
\\ &=
\omega \cdot \frac{N}{m \cdot (N-m)} \cdot \cov{ \att,\;  \sum_{i=1}^N Y_i(0) \cdot T_i  }
+
\omega \cdot \var{\att} 
\\ &\quad +
(1-\omega) \cdot \frac{N}{m \cdot (N-m)} \cdot \cov{ \atc ,\;  \sum_{i=1}^N Y_i(0) \cdot T_i  }  
+ 
(1-\omega) \cdot \cov{ \atc, \att } 
\end{align*}
According to the derivations that are described in detail later on in this section we can re-write the above covariance expression as:
\begin{align*}
\cov{ \ato, \; \td }
&=
\omega \cdot \frac{N}{m \cdot (N-m)} \cdot 
\frac{(1-p)N}{N-1} \cdot \left[ \rho \sigma_1 \sigma_0 - \sigma_0^2 \right]
+\omega \cdot \frac{N-m}{m(N-1)} \cdot \sigma_{\tau}^2 
\\ &\quad  
- (1-\omega) \cdot \frac{N}{m(N-m)} \cdot \frac{m}{N-1} \cdot \left[ \rho \sigma_1 \sigma_0 - \sigma_0^2 \right]
-(1-\omega) \cdot \frac{1}{N-1} \cdot \sigma_{\tau}^2
\\&=
\omega \cdot \frac{1}{p(N-1)} \cdot \left[ \rho \sigma_1 \sigma_0 - \sigma_0^2 \right]
+
\omega \cdot \frac{(1-p)}{p(N-1)} \cdot \sigma_{\tau}^2 - \frac{1}{N-1} (1-\omega) \sigma_{\tau}^2
\\ &\quad
-\frac{(1-\omega)}{ (N-1)(1-p) } \left[ \rho \sigma_1 \sigma_0 -\sigma_0^2 \right]
\\&=
\frac{\omega}{(N-1)p} \left[ \rho \sigma_1 \sigma_0 - \sigma_0^2 + (1-p) \sigma_{\tau}^2 + p \sigma_{\tau}^2 \right] - \frac{1}{N-1} \cdot \sigma_{\tau}^2
\\ &\quad
-\frac{(1-\omega)}{ (N-1)(1-p) } \left[ \rho \sigma_1 \sigma_0 -\sigma_0^2 \right]
\\&=
\frac{\omega}{(N-1)p} \left[ \sigma_1^2-\rho \sigma_1 \sigma_0   \right] - \frac{1}{N-1} \cdot \sigma_{\tau}^2
-\frac{(1-\omega)}{ (N-1)(1-p) } \left[ \rho \sigma_1 \sigma_0 -\sigma_0^2 \right]
\end{align*}

\subsection*{ Calculation of $\var{ \att}$ and $\var{ \atc}$ }

\begin{align*}
\var{ \att} &= \frac{N-m}{m(N-1)} \cdot \sigma_{\tau}^2 
=
\frac{N-m}{m(N-1)} \cdot \left[ \sigma_1^2 + \sigma_0^2 - 2 \rho \cdot \sigma_1 \cdot \sigma_0 \right]
\end{align*}
Similarly,
\begin{align*}
\var{ \atc} &= \frac{m}{(N-m)(N-1)} \cdot \sigma_{\tau}^2 
\end{align*}

\subsection*{ Calculation of $\cov{ \atc,\;  \sum_{i=1}^N Y_i(0) \cdot T_i  }$ }

\begin{align*}
\cov{ \atc,\;  \sum_{i=1}^N Y_i(0) \cdot T_i  } 
&= - \frac{m}{N-m} \cdot \cov{\att, \sum_{i=1}^N Y_i(0) \cdot T_i } 
\\ &= 
- \frac{m}{N-m} \cdot \frac{(1-p)N}{N-1} \cdot \left[ \rho \sigma_1 \sigma_0 - \sigma_0^2 \right] 
\\ &=
- \frac{m}{N-1} \cdot \left[ \rho \sigma_1 \sigma_0 - \sigma_0^2 \right] 
\\ &\approx
- p \cdot \left[ \rho \sigma_1 \sigma_0 - \sigma_0^2 \right] 
\end{align*}
where the second equality follows from the derivations in the next subsection below. 

\subsection*{ Calculation of $\cov{ \att,\;  \sum_{i=1}^N Y_i(0) \cdot T_i  }$ }
\begin{align*}
\cov{ \att,\;  \sum_{i=1}^N Y_i(0) \cdot T_i  } 
&= 
\frac{1}{m} \cdot \cov{ \sum_{i=1}^N Y_i(1) \cdot T_i - \sum_{i=1}^N Y_i(0) \cdot T_i ,\;  \sum_{i=1}^N Y_i(0) \cdot T_i  }
\\ &= 
\frac{1}{m} \cdot \left[ \cov{ \sum_{i=1}^N Y_i(1) \cdot T_i, \sum_{i=1}^N Y_i(0) \cdot T_i} 
- \var{ \sum_{i=1}^N Y_i(0) \cdot T_i } \right]
\\ &=
\frac{1}{m} \cdot \left[ \frac{p(1-p)N^2}{N-1} \cdot \rho \cdot \sigma_1 \sigma_0 
- \frac{m(N-m)}{N-1} \cdot \sigma_0^2 \right] 
\\ &=
\frac{(1-p)N}{N-1} \cdot \left[ \rho \sigma_1 \sigma_0 - \sigma_0^2 \right]
\end{align*}

\subsubsection*{ Derivations of $\cov{ \sum_{i=1}^N Y_i(1) \cdot T_i, \sum_{i=1}^N Y_i(0) \cdot T_i} $ }
\begin{align*}
\e{ \sum_{l=1}^N Y_l(1) \cdot T_l \cdot \sum_{i=1}^N Y_i(0) \cdot T_i  }
&= \e{ \sum_{i=1}^N \sum_{l=1}^N  Y_l(1) \cdot  Y_i(0) \cdot T_i \cdot T_l  }
\\ &=
\e{ T_i \cdot \sum_{i=1}^N Y_i(0) \cdot Y_i(1) + \sum_{i=1}^N \sum_{l \neq i} Y_i(0) \cdot Y_l(1) \cdot T_i \cdot T_l }
\\ &=
\Pr(T_i=1) \cdot \sum_{i=1}^N Y_i(0) \cdot Y_i(1) + \e{T_i \cdot T_l| i \neq l} \cdot \sum_{i=1}^N \sum_{l \neq i} Y_i(0) \cdot Y_l(1)
\\ &=
\frac{m}{N} \cdot \sum_{i=1}^N Y_i(0) \cdot Y_i(1) + \frac{m}{N} \cdot \frac{m-1}{N-1} \cdot \sum_{i=1}^N \sum_{l \neq i} Y_i(0) \cdot Y_l(1)
\\ &=
\frac{m}{N} \cdot \frac{N-m}{N-1} \cdot \sum_{i=1}^N Y_i(0) \cdot Y_i(1) + \frac{m}{N} \cdot \frac{m-1}{N-1} \cdot \sum_{i=1}^N \sum_{l=1 }^N Y_i(0) \cdot Y_l(1)
\end{align*} 
\begin{align*}
\e{ \sum_{l=1}^N Y_l(1) \cdot T_l } &= \frac{m}{N} \sum_{l=1}^N Y_l(1) \\
\e{ \sum_{l=1}^N Y_l(0) \cdot T_l } &= \frac{m}{N} \sum_{l=1}^N Y_l(0) \\
\Rightarrow 
\e{ \sum_{l=1}^N Y_l(1) \cdot T_l } \cdot \e{ \sum_{i=1}^N Y_l(0) \cdot T_i }
&=
\left( \frac{m}{N} \right)^2 \cdot \sum_{i=1}^N \sum_{l=1 }^N Y_i(0) \cdot Y_l(1)
\end{align*}
Hence,
\begin{align*}
\cov{ \sum_{i=1}^N Y_i(1) \cdot T_i, \sum_{i=1}^N Y_i(0) \cdot T_i} 
&= 
\frac{m}{N} \cdot \left[ \frac{N-m}{N-1} \cdot \sum_{i=1}^N Y_i(0) \cdot Y_i(1)  -\frac{N-m}{N(N-1)} \cdot \sum_{i=1}^N \sum_{l=1 }^N Y_i(0) \cdot Y_l(1)  \right]
\\ &= 
\frac{m}{N} \cdot  \frac{N-m}{N-1} \cdot  \left[ \sum_{i=1}^N Y_i(0) \cdot Y_i(1)  -\frac{1}{N} \cdot \sum_{i=1}^N \sum_{l=1 }^N Y_i(0) \cdot Y_l(1)  \right]
\\ &= 
\frac{p(1-p)N}{N-1} \cdot  \left[ \sum_{i=1}^N Y_i(0) \cdot Y_i(1)  -\frac{1}{N} \cdot \sum_{i=1}^N \sum_{l=1 }^N Y_i(0) \cdot Y_l(1)  \right]
\end{align*}
Re-arranging the above expression yields:
\begin{align*}
\cov{ \sum_{i=1}^N Y_i(1) \cdot T_i, \sum_{i=1}^N Y_i(0) \cdot T_i} 
&= \frac{p(1-p)N^2}{N-1} \cdot \rho \cdot \sigma_1 \sigma_0 
\end{align*}

\subsection*{ Calculation of $\cov{ \att,\;  \atc }$ }
Note that the covariance between SATT and SATC is,
\begin{align*}
\cov{\text{SATT}, \text{SATC} } &=
\cov{ \frac{1}{m} \cdot \sum_{i=1}^N (Y_i(1) - Y_i(0)) \cdot T_i,  \frac{1}{N-m} \cdot \sum_{i=1}^N (Y_i(1) - Y_i(0)) \cdot (1-T_i) }
\\ &= 
- \frac{1}{m} \cdot \frac{1}{N-m} \cdot \cov{ \sum_{i=1}^N (Y_i(1) - Y_i(0)) \cdot T_i,  \sum_{i=1}^N (Y_i(1) - Y_i(0)) T_i }
\\ &= 
- \frac{1}{N \cdot p (1-p)} \cdot \frac{\var{ \sum_{i=1}^N (Y_i(1) - Y_i(0)) T_i }}{N}
\\ &=
-\frac{1}{N \cdot p (1-p)} \cdot \frac{m^2}{N} \cdot \var{\att}
\\ &=
-\frac{1}{N \cdot p (1-p)} \cdot \frac{m^2}{N} \cdot \frac{N-m}{m(N-1)} \cdot \sigma_{\tau}^2
\\ &=
-\frac{1}{N-1} \cdot \sigma_{\tau}^2
\end{align*}

Another way of deriving the covariance between the SATT and SATC is to use the fact that $\var{\text{SATE}} = 0$:
First note that,
\begin{align*}
\var{\text{SATE}} &= \var{ \frac{m}{N} \cdot \text{SATT} } + \var{ \frac{N-m}{N} \cdot \text{SATC}} 
+ 2 \cdot p(1-p) \cdot \cov{\text{SATT}, \text{SATC}}
\\ &=
p^2 \cdot \frac{1}{m^2} \cdot \var{ \sum_{i=1}^N (Y_i(1) - Y_i(0)) T_i }
+
(1-p)^2 \cdot \frac{1}{(N-m)^2} \cdot \var{ \sum_{i=1}^N (Y_i(1) - Y_i(0)) T_i }
\\ & + 2 \cdot p(1-p) \cdot \cov{\text{SATT}, \text{SATC}} 
\\ &=
2 \cdot \frac{1}{N^2} \cdot \var{ \sum_{i=1}^N (Y_i(1) - Y_i(0)) T_i } 
+
2 \cdot p(1-p) \cdot \cov{\text{SATT}, \text{SATC}} 
\end{align*}
Hence,
\begin{align*}
\var{\text{SATE}} &= 0
\\ \Rightarrow 
\cov{\text{SATT}, \text{SATC}} &= - \frac{1}{N \cdot p (1-p)} \cdot \frac{\var{ \sum_{i=1}^N (Y_i(1) - Y_i(0)) T_i }}{N}
\end{align*}

\section{Proof of Theorem (\ref{theorem: CLT - Bernoulli trials randomization})}
\label{proof: CLT - Bernoulli trials randomization}

As both $m$ and $\sum_{i=1}^N Y_i(0) \cdot T_i$ are random variables, the distribution of the ratio can be derived using the Delta method. The result that the limiting distribution is Normal follows directly from standard results of the Delta method \citep{casella_berger}. The derivations of the variance-covariance matrix of the limiting distribution are detailed below. First notice that:
\begin{align*}
\frac{1}{N} \cdot \sum_{i=1}^N Y_i(0) \cdot T_i 
& \overset{p}{\rightarrow} 
\e{Y(0)} \\
\frac{1}{N} \cdot \sum_{i=1}^N Y_i(0) \cdot T_i 
& \overset{d}{\rightarrow} 
N \left( \e{Y(0)}, \;  \frac{1}{N} \cdot p(1-p) \left[ \sigma_0^2 + \e{Y(0)}^2 \right] \right)
\end{align*}
where 
\begin{align*}
\var{ \frac{1}{N} \cdot \sum_{i=1}^N Y_i(0) \cdot T_i } &=  \frac{1}{N^2} \cdot p(1-p) \sum_{i=1}^N Y_i(0)^2
\\ &=  
\frac{1}{N} \cdot p(1-p) \left[ \sigma_0^2 + \e{Y(0)}^2 \right]
\end{align*}
and that $m$ follows a Binomial distribution of $m \sim \text{Binom}(n, p)$ and the variance of $m$ is $\var{m} = N \cdot p(1-p)$. The covariance between $m$ and $ \sum_{i=1}^N Y_i(0) \cdot T_i $ is
\begin{align*}
\cov{m, \sum_{i=1}^N  Y_i(0) \cdot T_i } 
&= 
\cov{\sum_{i=1}^N  T_i , \sum_{i=1}^N  Y_i(0) \cdot T_i } 
= \sum_{i=1}^N  Y_i(0) \cdot \var{T_i} 
\nonumber \\ &=
p(1-p) \cdot \sum_{i=1}^N  Y_i(0) 
\nonumber \\ &=
p(1-p) \cdot N \cdot  \e{Y(0)}
\end{align*}   
and therefore the covariance-variance matrix is,
\begin{align}
\Sigma &= \left( 
\begin{array}{cc}
N \cdot p(1-p) \left[ \sigma_0^2 + \e{Y(0)}^2 \right] & p(1-p) \cdot N \cdot  \e{Y(0)} \\
p(1-p) \cdot N \cdot  \e{Y(0)} & N \cdot p(1-p) \\
\end{array}
\right)
\end{align}
The variance of $ \frac{1}{m} \cdot \sum_{i=1}^N Y_i(0) \cdot T_i $ is,
\begin{align*}
\left( \begin{array}{c}
\frac{1}{m} \\
 - \frac{\sum_{i=1}^N Y_i(0) \cdot T_i}{m^2}
\end{array} \right)^T
\cdot
\Sigma
\cdot
\left( \begin{array}{c}
\frac{1}{m} \\
 - \frac{\sum_{i=1}^N Y_i(0) \cdot T_i}{m^2}
\end{array} \right)
\end{align*}
where a consistent estimator for $p$ is $\frac{1}{N} \cdot \sum_{i=1}^N T_i$ and a consistent estimator for $\sum_{i=1}^N Y_i(0) \cdot T_i$ is $\sum_{i=1}^N Y_i(0) \cdot (1-T_i)$.

\end{document}